\documentclass[12pt,reqno,twoside]{amsart}

\usepackage{amsfonts}
\usepackage{amsmath,amssymb,amsthm,amsthm}
\usepackage{mathtools}
\usepackage{dsfont}
\usepackage{etoolbox}
\usepackage{color}
\usepackage{bm}
\usepackage{units}
\usepackage{enumitem}
\usepackage[skip=0.333\baselineskip]{subcaption}
\usepackage{caption}

\usepackage{hyperref,graphicx}
\hypersetup{%
        colorlinks,
        breaklinks=true,
        plainpages=false,
        citecolor=.,
        linkcolor=.,
        urlcolor=.,
        bookmarksopen=true,%
        bookmarksnumbered=true,
        bookmarksdepth=8,
        unicode=true,%
        pdfpagelayout={SinglePage},pdffitwindow=true}%

\usepackage{cleveref}
\usepackage{booktabs}

\usepackage[dvips,bottom=1.4in,right=1in,top=1in, left=1in]{geometry}

\numberwithin{figure}{section}
\numberwithin{table}{section}

\makeatletter
\def\eqnarray{\stepcounter{equation}\let\@currentlabel=\theequation
\global\@eqnswtrue
\tabskip\@centering\let\\=\@eqncr
$$\halign to \displaywidth\bgroup\hfil\global\@eqcnt\z@
  $\displaystyle\tabskip\z@{##}$&\global\@eqcnt\@ne
  \hfil$\displaystyle{{}##{}}$\hfil
  &\global\@eqcnt\tw@ $\displaystyle{##}$\hfil
  \tabskip\@centering&\llap{##}\tabskip\z@\cr}

\def\endeqnarray{\@@eqncr\egroup
      \global\advance\c@equation\m@ne$$\global\@ignoretrue}

\usepackage{placeins}

\let\Oldsection\section
\renewcommand{\section}{\FloatBarrier\Oldsection}

\def\xvec{{\bf x}}
\def\fvec{{\bf f}}
\def\uvec{{\bf u}}
\def\svec{{\bf s}}

\def\Imat{{\bf I}}
\def\Kmat{{\bf K}}

 \usepackage[textsize=small]{todonotes}
\title[Risk Measures in Digital Twins]{On The Use of Risk Measures in Digital Twins to \\ Identify Weaknesses in Structures}

\author[Airaudo]{Facundo N. Airaudo}
\address{Center for Computational Fluid Dynamics and 
Department of Physics and Astronomy, 
George Mason University, Fairfax, VA 22030, USA.}
\email{fairaudo@gmu.edu}
\author[Antil]{Harbir Antil}
\address{Center for Mathematics and Artificial Intelligence (CMAI) and 
Department of Mathematical Sciences, 
George Mason University, Fairfax, VA 22030, USA.}
\email{hantil@gmu.edu}
\author[L\"ohner]{Rainald L\"ohner}
\address{Center for Computational Fluid Dynamics and 
Department of Physics and Astronomy, 
George Mason University, Fairfax, VA 22030, USA.}
\email{rlohner@gmu.edu}
\author[Rakhimov]{Umarkhon Rakhimov}
\address{Center for Mathematics and Artificial Intelligence (CMAI) and 
Department of Mathematical Sciences, 
George Mason University, Fairfax, VA 22030, USA.}
\email{urakhimov@gmu.edu}

\begin{document}

\maketitle

\begin{abstract}
Given measurements from sensors and a set of standard forces, 
an optimization based approach to identify weakness 
in structures is introduced. The key novelty lies in letting the load and 
measurements to be random variables. Subsequently the 
conditional-value-at-risk (CVaR) is minimized subject to the elasticity 
equations as constraints. CVaR is a risk measure that leads to minimization
of rare and low probability events which the standard expectation cannot. 
The optimization variable is the (deterministic) strength factor which 
appears as a coefficient in the elasticity equation, thus making the problem
nonconvex. Due to uncertainty, the problem is high dimensional and, due 
to CVaR, the problem is nonsmooth. An adjoint based 
approach is developed with quadrature in the random variables. Numerical 
results are presented in the context of a plate, a large structure with trusses similar to those used in solar arrays or cranes, and a footbridge.
\end{abstract}

\section{Introduction}

Given that all materials exposed to the environment and/or undergoing
loads eventually age and fail, the task of trying to detect and localize
weaknesses in structures is common to many fields. To mention just a few:
airplanes, drones and missiles, turbines, launch pads and airport
infrastructure, wind turbines, satellites and space stations.
Traditionally, manual inspection was the only way of carrying out this task,
aided by ultrasound, X-ray, or vibration analysis techniques.
The advent of accurate, abundant and cheap sensors, together with
detailed, high-fidelity computational models
in an environment of digital twins has opened the possibility of
enhancing and automating the detection and localization of
weaknesses in structures. The authors in article \cite{FAiraudo_RLoehner_HAntil_2023a} introduced an adjoint based optimization approach \cite{JLLions_1971a,FTroeltzsch_2010a,MHinze_RPinnau_MUlbrich_SUlbrich_2009a,HAntil_DPKouri_MDLacasse_DRidzal_2018a} to identify weakness in structures under the deterministic setting.
Given the displacement or strain measurements 
from certain sensor measurements, the goal is to solve an inverse problem to 
determine material properties. This amounts to minimizing a cost functional 
subject to elasticity equations as constraints, and is a deterministic optimization 
problem with partial differential equations (PDE) as constraints. 
Here the optimization problem itself can be thought as a `digital twin' which is informed by the physical system via sensor measurements. The digital twin then make predictions about the structural weakness.

However, the underlying elasticity equation contains various unknown quantities. In particular, the load measurements can present uncertainty, especially if it comes from sources that were not accounted for, like wind or temperature variations. In addition, the measurements from sensors
could be faulty due to sensor errors or various signal-to-noise thresholds. 
The present article proposes to model 
these quantities as random fields. This has dual benefits: firstly, one can tackle the 
unknown quantities and secondly, it will lead to designs which are resilient to 
uncertainty. However, several challenges appear. Due to random data/inputs, the
PDE solution becomes a random field and the cost functional becomes a random
variable. It is unrealistic to minimize a random variable. Moreover, these problems can be very high dimensional problems. 

Traditionally, cost functionals that are random variables have been handled 
by computing their expectation. However, evaluating the expectation is like computing 
an average. This scenario cannot handle the outliers which lie in the tail of a 
distribution. Instead, following \cite{AShapiro_DDentcheva_ARuszczynski_2014a}, we consider 
the so-called conditional-value-at-risk (CVaR), which is the average of the $\beta$-tail 
of the distribution with $\beta \in (0,1)$. The first use of CVaR in optimization problems 
with PDE constraints can be found in \cite{DPKouri_TMSurowiec_2016a}. Recently, 
this approach has been further extended using tensor decomposition tools 
\cite{HAntil_SDolgov_AOnwunta_2022b} to tackle higher dimensional problems. 
See also \cite{HAntil_SDolgov_AOnwunta_2023a} for a tensor-train framework with state
constraints (e.g., constraints on the displacement).

\medskip
\noindent
{\bf Outline:} The remainder of the paper has been organized as follows: In 
section~\ref{s:cvar} we state the well known definition of CVaR. 
Section~\ref{s:problem} focuses on the optimization problem with elasticity
constraints under uncertainty. Section~\ref{s:adj} discusses the adjoint formulation to evaluate the gradients. Section~\ref{s:numerics}, presents our implementation details and various numerical examples. Finally, in section~\ref{s:conclusion} a few concluding remarks are provided.

\section{Conditional Value at Risk}
\label{s:cvar}

Let $\mathcal{L} := (D,\mathcal{A},\mathbb{P})$ be a complete probability space.
Here $D$ denotes the set of outcomes, $\mathcal{A} \subset 2^D$ is the $\sigma$-algebra
of events, and $\mathbb{P} : \mathcal{A} \rightarrow [0,1]$ is the appropriate 
probability measure. Consider a scalar random variable $X$ defined on $\mathcal{L}$. In our setting below, $X$ will be the random variable objective function.
The expectation of $X$ is given by
\begin{equation}\label{eq:risk_neutral}
		\mathbb{E}[X] = \int_D X(\omega) d\mathbb{P}(\omega) \, .
\end{equation}
$X$ can be, for example, a typical objective function in an optimization problem, minimizing its standard expectation would be regarded as a risk neutral approach. Making use of risk measures like CVaR$_\beta$ allow the optimization process to focus on certain areas of the random domain instead of the whole range.

\begin{figure}[!hbt]
    \centering
    \includegraphics[width=0.7\textwidth]{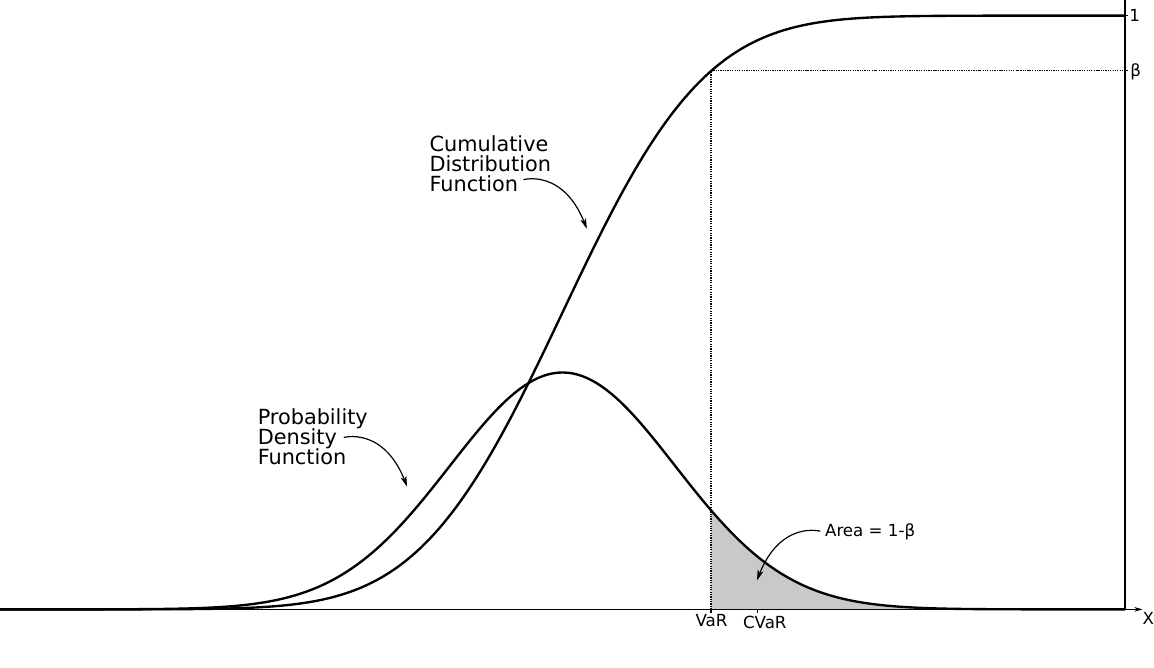}
    \caption{Illustrations of a normal Probability Density Function (PDF) and its corresponding Cumulative Distribution function (CDF) for a random loss function $X$. For an arbitrary $\beta$,  VaR$_\beta$ and CVaR$_\beta$ are shown.}
    \label{fig:cvar_plot}
\end{figure}

If $\beta \in (0,1)$ is fixed, and $\mathbb{P}[X \le t] $ denotes the probability of the 
random variable $X$ is less than or equal to $t$, then the value-at-risk (VaR) 
is given by 
\begin{equation} \label{eq:VaR}
		\mbox{VaR}_\beta[X] := \inf_{t\in \mathbb{R}} \{ t \, : \, \mathbb{P}[X \le t] \ge \beta \} \, .
\end{equation}
Unfortunately, $\mbox{VaR}_\beta$ is not a coherent risk measure as it violates the 
sub-additivity / convexity axiom in risk-measures \cite{AShapiro_DDentcheva_ARuszczynski_2014a}.  
Instead, $\mbox{CVaR}_\beta$, which is coherent, is a preferred risk-measure. 
Due to Rockafellar and Uryasev \cite{TRRockafellar_SUryasev_2000a,TRRockafellar_SUryasev_2002a}, 
$\mbox{CVaR}_\beta$ can be written as 
	\begin{equation}\label{eq:CVaR}
		\mbox{CVaR}_\beta[X] = \inf_{t\in \mathbb{R}} 
						\left\{  t + \frac{1}{1-\beta} \mathbb{E}[(X-t)_+] \right\} \, ,
	\end{equation}
where $(X-t)_+ = \max\{ X - t , 0 \}$. In case of a continuous random variable, the
above definition is equivalent to $\mbox{CVaR}_\beta[X] = \mathbb{E}[X \, : \, X > \mbox{VaR}_\beta[X]]$,
i.e., $\mbox{CVaR}_\beta[X]$ is indeed the average of $\beta$-tail of the distribution
of $X$. Namely, $\mbox{CVaR}_\beta[X] $ focuses on the rare and low probability
events, especially when $\beta \rightarrow 1$. This is further illustrated in Figure~\ref{fig:cvar_plot} by the shaded region under the curve. I.e., the Cumulative Distribution Function (CDF) here is larger than $\beta$ and here $X > \mbox{VaR}_\beta[X]$ with probability $1-\beta$. Notice that the main reason for considering $\mbox{CVaR}_\beta$ instead of $\mbox{VaR}_\beta$ is that the latter is not coherent. In addition, the expression of $\mbox{CVaR}_\beta$ given in \eqref{eq:CVaR} making it much more computationally tractable in comparison to $\mbox{VaR}_\beta$. 

Throughout, we make the finite dimensional noise assumption. We assume that 
$\omega$ can be sampled via a finite random vector $\xi : D \rightarrow \Xi$ instead,
where $\Xi = \xi(D) \subset \mathbb{R}^d$ with $d \in \mathbb{N}$. This allows us
to redefine the probably space as $(\Xi,\Sigma,\gamma)$ where $\Sigma = \xi(\mathcal{A})$
is the $\sigma$-algebra and $\gamma(\xi)$ is the continuous probability density function
such that $\mathbb{E}[X] = \int_\Xi X(\xi) \gamma(\xi) d\xi$. The random variable $X(\xi)$ can
be considered as a function of the random vector $\xi = (\xi^{(1)}, \dots, \xi^{(d)})$. 







\section{Problem Formulation}
\label{s:problem}

In what follows we assume that a spatial discretization has been carried out
via finite elements, finite volumes or finite differences, even though the discussion below is independent of any particular discretization.
The determination of material properties (or weaknesses)
may be formulated as an optimization problem for the deterministic 
strength factor $\alpha(\xvec)$ as follows: 
Let the load $\Xi \ni \xi \mapsto \fvec(\xi)$ be a random vector with known 
distribution. Moreover, let $\xvec_j, j=1,\dots,m$ denotes the $m$ 
measurement locations for deformations $\Xi \ni \xi \mapsto \uvec^{md}_{j}(\xi), \ j=1,\dots,m$ 
or strains $\Xi \ni \xi \mapsto \svec^{ms}_{j}(\xi), \ j=1,\dots,m$. 

Then the random variable objective functional that we want to minimize
to identify the spatial distribution of the strength factor $\alpha$ is given
by 
\begin{equation} \label{eq:cost0}
 I(\uvec(\xi),\alpha) = 
  {1 \over 2}  \sum_{j=1}^m w^{md}_{j}(\xi) 
             ( \uvec^{md}_{j}(\xi) - \Imat^d_{j} \cdot \uvec(\xi) )^2 
+ {1 \over 2}  \sum_{j=1}^m w^{ms}_{j}(\xi) 
             ( \svec^{ms}_{j}(\xi) - \Imat^s_{j} \cdot \svec(\xi) )^2 \, ,
\end{equation}
where $w^{md}_{j}, w^{ms}_{j}$ are displacement and strain weights
possibly depending on parameter $\xi$, 
$\Imat^d, \Imat^s$ interpolation matrices that are used to obtain 
the displacements and strains from the finite element mesh
at the measurement locations. 
However, as mentioned earlier, minimizing $I(\uvec(\xi),\alpha) $ in \eqref{eq:cost0} 
is not tractable. 
We scalarize this using CVaR$_\beta$ given in \eqref{eq:CVaR}. The resulting 
minimization problem is given by 
\begin{subequations}
\begin{equation} \label{eq:cost0}
 	\inf_{\alpha, t} \left\{ \mbox{CVaR}_\beta[I(\uvec(\cdot),\alpha)] 	
		:= \left\{  t + \frac{1}{1-\beta} \mathbb{E}[(I(\uvec(\cdot),\alpha)-t)_+] \right\} \right\}
\end{equation}
subject to the finite element description (e.g. trusses,
beams, plates, shells, solids) of the structure \cite{zienkiewicz2005finite, simo2006computational} under 
consideration (i.e. the digital twin/system \cite{mainini2015surrogate, chinesta2020virtual}):
\begin{equation} \label{eq:forward_pde}
    \Kmat(\alpha) \uvec(\xi) = \fvec(\xi) ,  \; \quad \mbox{a.s. } \xi \in \Xi \, .
\end{equation}
\end{subequations}
Here $\Kmat$ is the usual stiffness matrix, which is obtained by assembling 
all the element matrices:
\begin{equation} \label{eq:strength_factor}
    \Kmat = \sum_{e=1}^{N_e} \alpha_{e} \Kmat_{e} \, ,
\end{equation}
where the strength factor $\alpha_{e}$ of the elements has already been
incorporated. We note in passing that in order to ensure that $\Kmat$
is invertible and non-degenerate $\alpha_{e} > \epsilon > 0$.

\section{Adjoint Approach}
\label{s:adj}

In order to establish a gradient based method, we need to calculate 
the derivative of CVaR$_\beta$ with respect to $\alpha$. Notice that 
$(\cdot)_+$ is not differentiable in the classical sense, but it is 
differentiable in a generalized sense \cite{MUlbrich_2011a,HAntil_LBetz_DWachsmuth_2023a}.  
We set the derivative of $(x)_+$ as 
\[
	(x)_+' = 
	\left\{
	\begin{array}{cc}
		1  & \mbox{if } x \ge 0 \\
		0  & \mbox{otherwise} \, ,
	\end{array}
	\right.
\]
enabling us to use a gradient based method. 
If one wants to use a higher order method (like Newton), 
then one approach is to smooth the $(\cdot)_+$ function 
\cite{DPKouri_TMSurowiec_2016a,HAntil_SDolgov_AOnwunta_2023a}.

Next, let the Lagrangian functional be  
\[
	\mathcal{L}(\uvec,\alpha,\widetilde{\uvec}) 
	= \mbox{CVaR}_\beta[I(\uvec(\cdot),\alpha)] 	
	   + \int_\Xi \gamma(\xi) \, \widetilde\uvec(\xi)^\top ( \Kmat(\alpha) \uvec(\xi) - \fvec(\xi) ) d\xi 
\]
where $\widetilde\uvec$ indicates the Lagrange multiplier. A variation of
$\mathcal{L}$ with respect to $\uvec$, at a stationary point, leads to 
\[
	0 = \delta_\uvec \mathcal{L}(\uvec,\alpha,\widetilde{\uvec}) 
	   = \frac{1}{1-\beta} \mathbb{E}[ (I(\uvec(\cdot),\alpha)-t)_+' \delta_\uvec I(\uvec(\cdot),\alpha) ]
	       +   \int_\Xi \gamma(\xi) \, \Kmat(\alpha)^\top \widetilde\uvec(\xi) d\xi \, ,
\]
which gives rise to the adjoint equation
\begin{equation}\label{eq:adj_eqn}
	\Kmat(\alpha)^\top \widetilde\uvec(\xi) = \frac{-1}{1-\beta} 
						\left[ (I(\uvec(\xi),\alpha)-t)_+' \ \delta_\uvec I(\uvec(\xi),\alpha) \right] 
                ,  \; \quad \mbox{a.s. } \xi \in \Xi \, .
\end{equation}
Subsequently, the variation of $\mathcal{L}$ with respect to $\alpha$ 
gives us the required derivative with respect to $\alpha$
\begin{equation}\label{eq:DaL}
	\delta_\alpha \mathcal{L}(\uvec,\alpha,\widetilde{\uvec}) 
	= \int_\Xi \gamma(\xi) \, \widetilde\uvec(\xi)^\top \Kmat'(\alpha) \uvec(\xi) d\xi .
\end{equation}
Finally, the variation with respect to $t$ is given by 
\begin{equation}\label{eq:Dat}
	\delta_t \mathcal{L}(\uvec,\alpha,\widetilde{\uvec}) 
	= 1 - \frac{1}{1-\beta} \mathbb{E} \left[(I(\uvec(\cdot),\alpha)-t)_+' \right] \, .
\end{equation}
With the expression of the derivatives of $\mathcal{L}$ with respect to 
$\alpha$ and $t$, we can now create a gradient based method to solve
the optimization problem.

\section{Implementation Details and Numerical Experiments}
\label{s:numerics}

The goal of this section is provide several illustrative experiments. 
Prior to that, we discuss some missing ingredients. Section~\ref{s:weight}
discusses the appropriate weights to be used in \eqref{eq:cost0}. 
This is followed by section~\ref{s:smoothing} which discusses 
the appropriate smoothing procedure for the gradients, first 
introduced in \cite{FAiraudo_RLoehner_HAntil_2023a}. The article 
\cite{FAiraudo_RLoehner_HAntil_2023a} also highlights that this 
smoothing can be associated with the use of appropriate function  
space based scalar product. 
Section~\ref{s:Eapprox}
focuses on the approximation of the expectation $\mathbb{E}$.
The proposed approach is applied to illustrative examples
in section~\ref{s:plate}. There examples were computed using our 
codes for optimization and interfacing, and CALCULIX 
\cite{dhondt2022calculix} as our solver for the elasticity equations. CALCULIX is a general, open source 
finite element code for structural mechanics applications with 
many element types, material models and options.

\subsection{Local Weighting}
\label{s:weight}

In this study, we consider the following weights in \eqref{eq:cost0},
and refer to \cite{FAiraudo_RLoehner_HAntil_2023a} for further 
scenarios: 
\begin{equation} \label{eq:local_weighting}
    w^{md}_{ij}={1 \over{(\uvec^{md}_{ij})^2}} ~~;~~ 
    w^{ms}_{ij}={1 \over{(\svec^{ms}_{ij})^2}} ~~;
\end{equation}

\subsection{Smoothing of Gradients} \label{s:smoothing}

The gradients of the cost function with respect to $\alpha$ allow for
oscillatory solutions. One must therefore smooth or `regularize' the
spatial distribution. This happens naturally when using few degrees of
freedom, i.e. when $\alpha$ is defined via other spatial shape functions
(e.g. larger spatial regions of piecewise constant $\alpha$). As the
(possibly oscillatory) gradients obtained in the (many) finite elements
are averaged over spatial regions, an intrinsic smoothing occurs.
This is not the case if $\alpha$ and the gradient are defined and 
evaluated in each element separately, allowing for the largest degrees of
freedom in a mesh and hence the most accurate representation.
Several types of smoothing or `regularization' are possible, see 
\cite{FAiraudo_RLoehner_HAntil_2023a}. 
All of them start by performing a volume averaging from elements
to points:
\begin{equation} \label{eq:averaging}
    \alpha_p = {{ \sum_{e} \alpha_{e} V_{e} } \over
               { \sum_{e} V_{e} }}
\end{equation}
where $\alpha_p, \alpha_{e}, V_{e} $ denote the value of $\alpha$ at
point $p$, as well as the values of $\alpha$ in element $e$ and the
volume of element $e$, and the sum extends over all the elements
surrounding point $p$. This work uses the averaging process 
described next.

\subsubsection*{Simple Point/Element/Point Averaging}

In this case, the values of $\alpha$ are cycled between elements 
and points. When going from point values to element values, 
a simple average is taken:
\begin{equation} \label{eq:simple_averaging}
    \alpha_{e} = { 1 \over n_{e} } \sum_i \alpha_i
\end{equation}
where $n_{e}$ denotes the number of nodes of
an element and the sum extends over all the nodes of the element.
After obtaining the new element values via equation \eqref{eq:simple_averaging} the point
averages are again evaluated via equation \eqref{eq:averaging}. These two steps are
repeated for a number of smoothing steps, which typically are in the range 1-5.
This form of averaging is very crude, but works very well in practice.

\subsection{Approximation of expectation}
\label{s:Eapprox}

The expectation $\mathbb{E}[\cdot]$ in \eqref{eq:cost0}, \eqref{eq:DaL}, 
and \eqref{eq:Dat} is approximated using Gauss quadrature. Unfortunately, 
the computational complexity in this case is $\mathcal{O}(n^d)$ where $n$ 
is the number of quadrature points in each random variable direction and 
$d$ is dimension of the random variables. One can overcome this curse of 
dimensionality by using tensor-train (TT) decomposition as recently shown 
in \cite{HAntil_SDolgov_AOnwunta_2023a,HAntil_SDolgov_AOnwunta_2022b}. 
However, the TT decomposition is not used in the present paper.

\begin{figure}[!hbt]
    \centering
    \includegraphics[width=\textwidth]{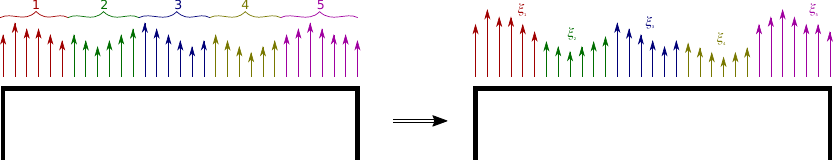}
    \caption{Illustration of the use of load groups. Each group is scaled by a different random variable $\xi_i$.}
    \label{fig:load_groups}
\end{figure}
A simpler approach is chosen, in the form of load groups. Since we expect the load to be distributed in space along a certain direction, we divide the load factor into multiple groups, each scaled separately by a random factor. Figure \ref{fig:load_groups} shows an example of this concept. The choice of how many groups to use is given by a compromise between required precision and availability of computing resources. In our computations, we choose the groups of uniform width.

For the cases in which a certain degree of uncertainty is expected in the measured loads, 
we have the following formulation for the forward problem
\begin{equation} \label{eq:forward_pde_load_groups}
    \Kmat(\alpha) \uvec = \mathbf{Q}_\xi \fvec,
\end{equation}
where $\mathbf{Q}_\xi$ is a diagonal matrix that makes it so the load applied to point $i$ is scaled by the proper random variable $\xi_{i_g}$.

For example, the standard expectation for a continuous random variable from equation \eqref{eq:risk_neutral}, rewritten here for clarity:
\begin{equation*}
    \mathbb{E}[X] = \int_D X(\omega) d\mathbb{P}(\omega).
\end{equation*}

Recalling the finite dimensional noise assumption from section \ref{s:cvar}, we can rewrite it as
\begin{equation*}
    \mathbb{E}[X] = \int_\Xi X(\xi) \gamma(\xi) d\xi,
\end{equation*}
where $\xi$ is a vector $\xi = (\xi_1, ..., \xi_{n_l})$ representing the load groups with
$n_l$ indicating the number of load groups.

If we discretize the integral using a Gauss quadrature, we have
\begin{equation}
    \mathbb{E}[X] = \sum_{i_g^{(1)}=1}^{n_g} ... \sum_{i_g^{(n_l)}=1}^{n_g} w_{i_g^{(1)}} ... w_{i_g^{(n_l)}} X(x_{i_g^{(1)}},...,x_{i_g^{(n_l)}}) \gamma(x_{i_g^{(1)}},...,x_{i_g^{(n_l)}}).
\end{equation}
For a quadrature with $n_g$ points in each direction, and using $n_l$ load groups, we need to compute $X$ and $\gamma$ a total of $n_g^{n_l}$ times. Clearly, this can quickly become very costly.

\subsection{Example: Plate.}
\label{s:plate}
\begin{figure}[!hbt]
    \centering
    \includegraphics[width=0.46\textwidth]{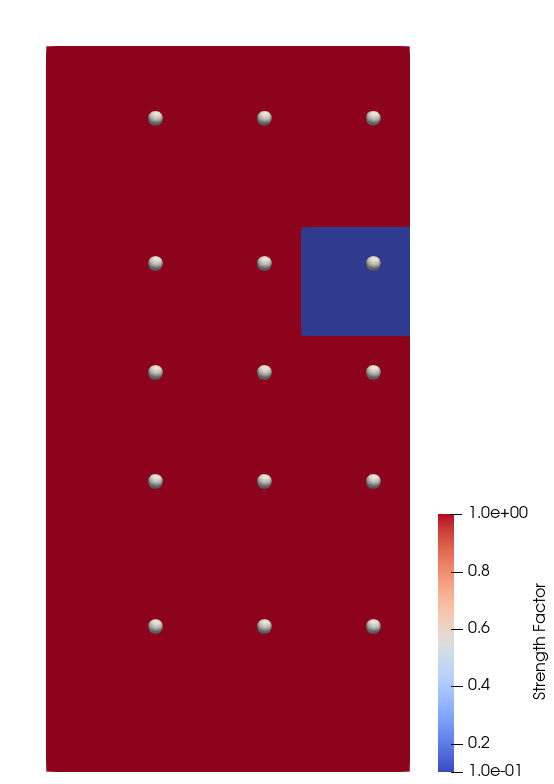} \quad
    \includegraphics[width=0.46\textwidth]{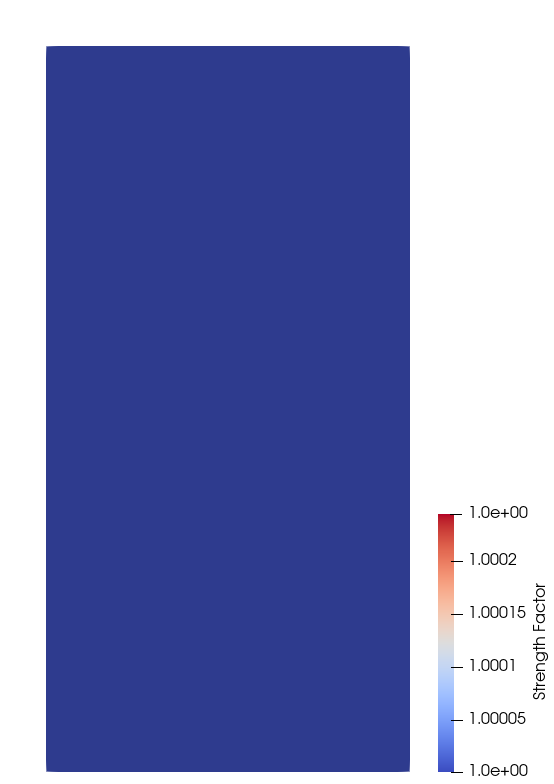}
    \caption{Left: Target strength factor and sensors. Right: Initial Strength Factor.}
    \label{fig:plate_target_initial}
\end{figure}

%
We consider the plate shown in Figure~\ref{fig:plate_target_initial}. The plate dimensions 
are (all units in mks): $0 \le x \le 60$, $0 \le y \le 30$, $0 \le z \le 0.1$.
Density, Young’s modulus and Poisson rate were set to $\rho = 7.8$, 
$E = 2 \cdot 10^{9} $, $\nu = 0.3$, respectively. The top boundary is clamped
$(\uvec = 0)$ and a vertical load is given by $f_y = 4 \cdot 10^5 \xi$, where the scalar $\xi$ 
lies in $\Xi = (0.8,1.2)$.

The spatial discretization for the displacement and adjoint is carried out 
using piecewise linear finite elements. The strength factor is discretized
using piecewise constant elements. 

Gauss quadrature with 4 terms is used to approximate the integrals 
over $\Xi$. Steepest descent with backtracking line search is used as the 
optimization algorithm. Finally, four smoothing steps are applied to the 
gradient.

Consider the configuration shown in Figure~\ref{fig:plate_target_initial} 
where part of plate has
been weakened. The goal is to try to match the displacement corresponding to
the configuration shown in this figure. 

\subsubsection{With Load and Sensor Uncertainty}

For each $j$, the random  
load is $\hat{f}_y = 10^5 \xi$ where $\xi$ is randomly drawn from a uniform 
distribution over $\Xi$. For each $j$, this gives us the desired measurements 
$\uvec_j^{md}$. That means that every sensor will derive its values from a different load case. Thus, we can see it as having $m$ conflicting realizations of the random variable at the same time.

\begin{figure}[h!]
        \centering
    \includegraphics[width=0.25\textwidth]{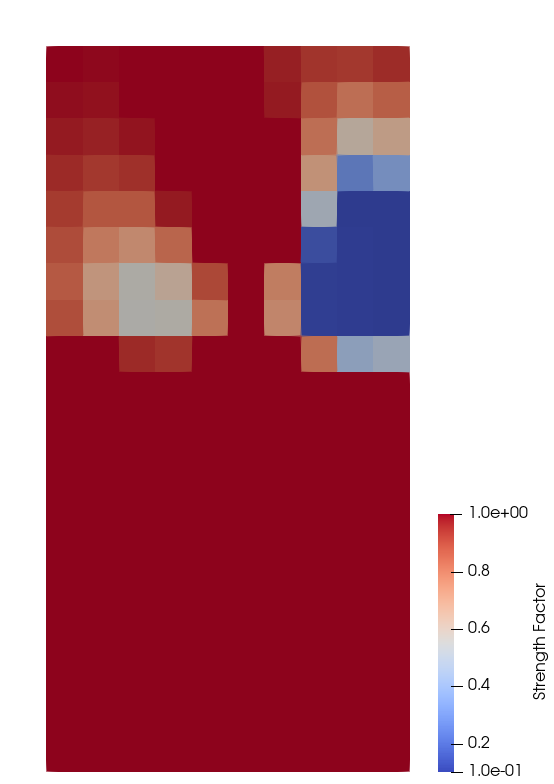}  \qquad      
	\includegraphics[width=0.25\textwidth]{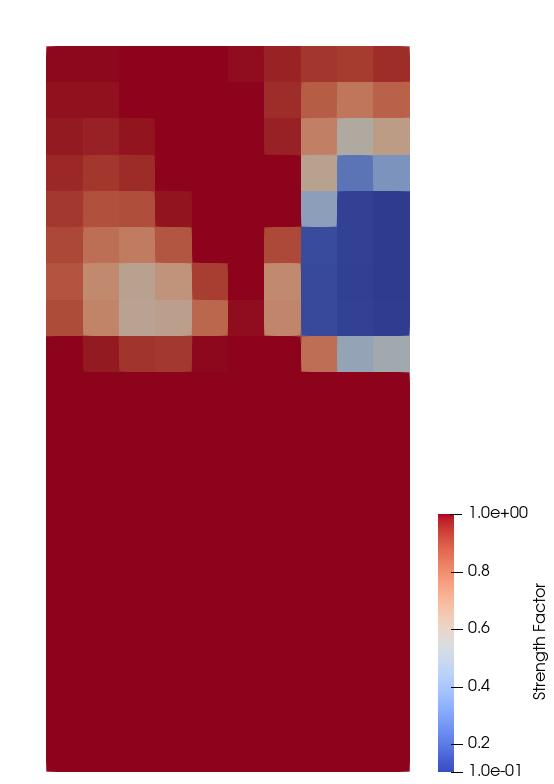}\qquad
	\includegraphics[width=0.25\textwidth]{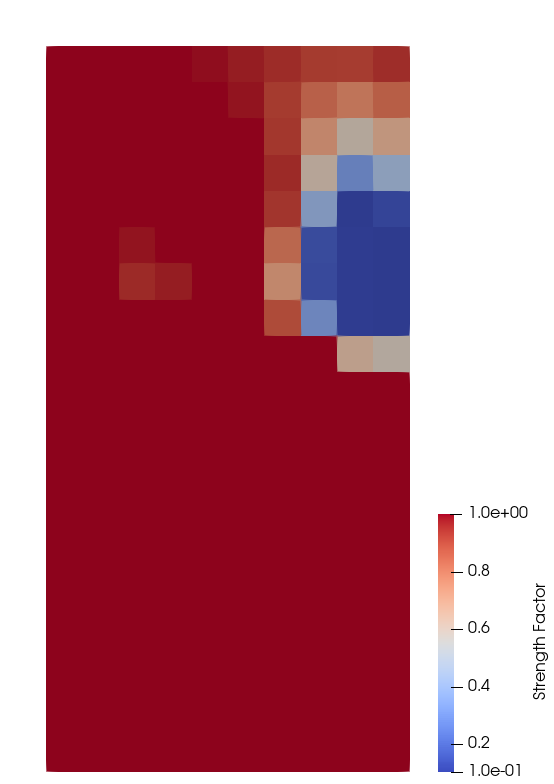}\qquad
	\includegraphics[width=0.25\textwidth]{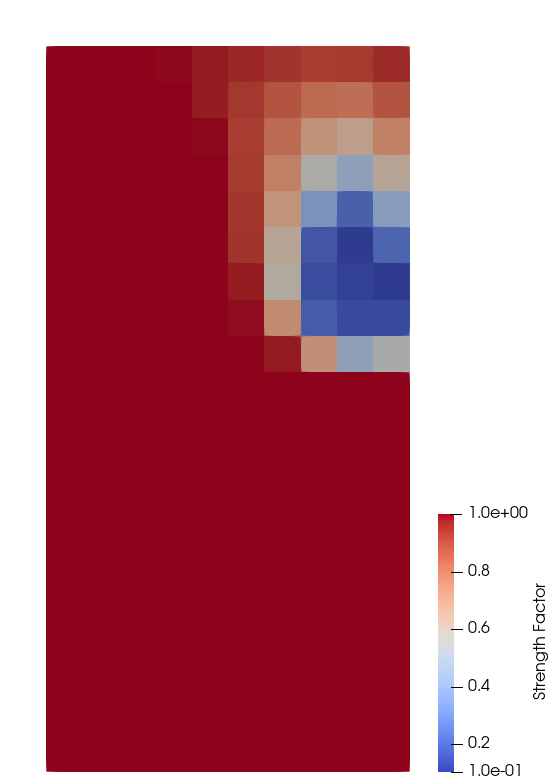}\qquad
	\includegraphics[width=0.25\textwidth]{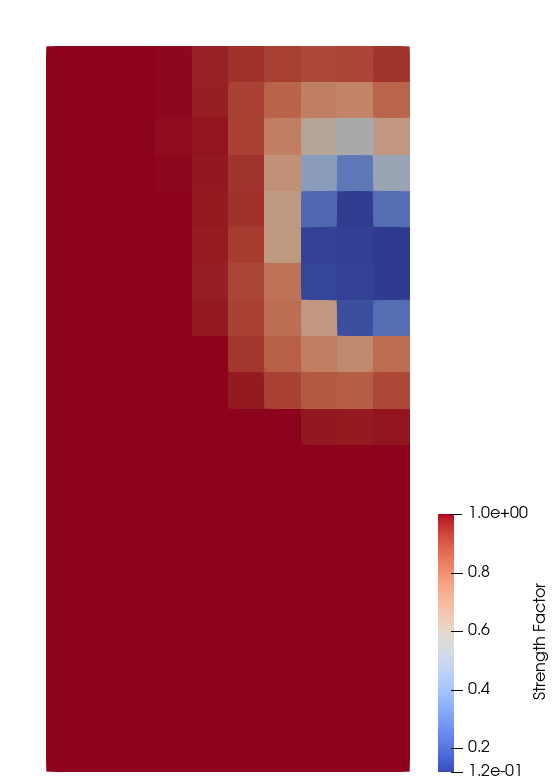}\qquad
    \includegraphics[width=0.25\textwidth]{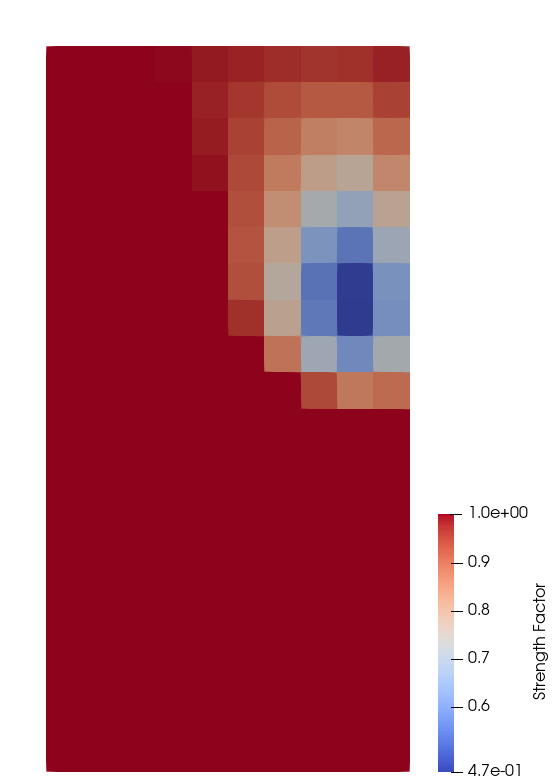}
	\caption{
    Top row: first panel (standard expectation); middle and right panels (CVaR$_\beta$ with
    $\beta = 0.1$ and $0.3$). Bottom row (from left to right): CVaR$_\beta$ with 
    $\beta = 0.5$, $0.7$ and $0.9$. 
    Clearly, CVaR$_\beta$ performs much better than the standard expectation.}
        \label{fig:plate_results}
\end{figure}

Figure~\ref{fig:plate_results} shows a comparison 
between the standard expectation and CVaR$_\beta$ with different levels of $\beta$. All optimization runs took 100~iterations for a decrease in the objective function of 2~orders of magnitude. Clearly, CVaR$_\beta$ performs much better as it is able to get rid of the extra weak spot that the uncertainties produced.

It is worth noting that as larger $\beta$ are utilized, the solution tends to become more underdeveloped, meaning that the values of the parameter do not go all the way down to $0.1$ as we would like. However, this is not a problem since our main goal is to accurately find the location of the weakening and the value placed on it is not as important.

\subsubsection{With Load Uncertainty}

In this example, we look at high dimensionality in the load uncertainty. For this case, all sensors will 
measure the same (deterministic) constant load given by $f_y = 4 \cdot 10^5$, applied at the top surface. In 
our optimization problem, we apply a random load as described below.

\begin{figure}[!hbt]
    \centering
    \includegraphics[width=0.95\textwidth]{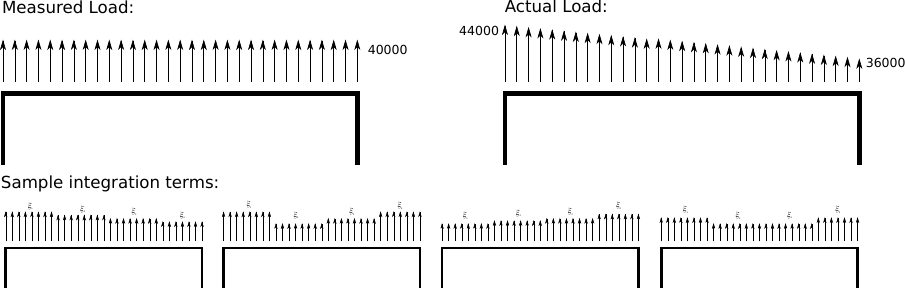}
    \caption{Top Left: Loading assumed to be present. Top Right: Actual loading used to compute the sensor displacements. Bottom: Examples of loading configurations used in the integration process.}
    \label{fig:plate_random_setup}
\end{figure}

In order to properly test a high-dimensional, worst case scenario, we compute the target displacements using a deterministic load with a linear variation $f_y = 4.4 \cdot 10^5 - \frac{8 \cdot 10^4}{60} x$, as shown in Figure \ref{fig:plate_random_setup} (top right). This sets up an interesting dimensionality problem, as the load was assumed to be constant, the linear behaviour of this load is not known by the optimizer.

The best one can do if one expects a non-constant variation of the loading in space, is to set up as many load groups as one can afford to in order to best adapt to unforeseen random behaviour. In this example, we choose to divide the loaded surface in 4 groups, each with its own associated random variable. Illustrations of this concept can be seen in the bottom row of figure \ref{fig:plate_random_setup}.



Using a Gauss quadrature scheme of order 3, our optimization algorithm required $3^4$ state solves in order to evaluate the objective function and the same number of adjoint solves to evaluate the gradient during each optimization iteration. Recall that the quadrature is used to approximate the expectation in both the objective function and the gradient with respect to $\alpha$.

\begin{figure}[!hbt]
    \centering
    \includegraphics[width=0.25\textwidth]{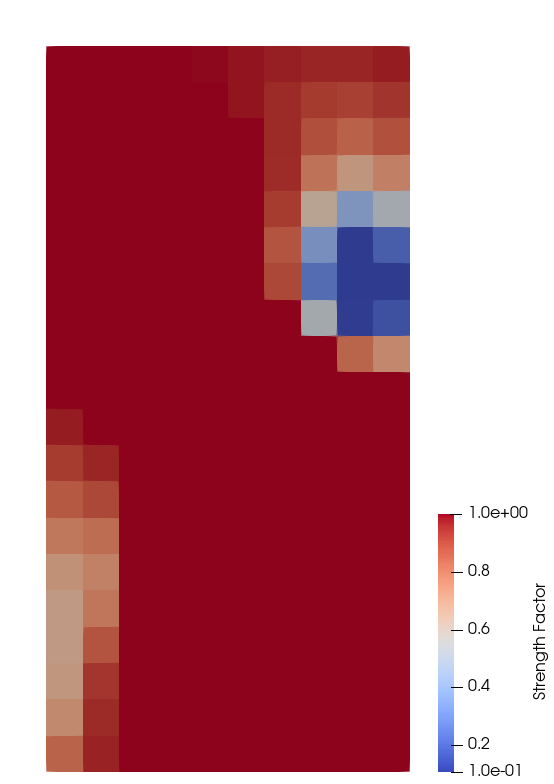}  \qquad      
	\includegraphics[width=0.25\textwidth]{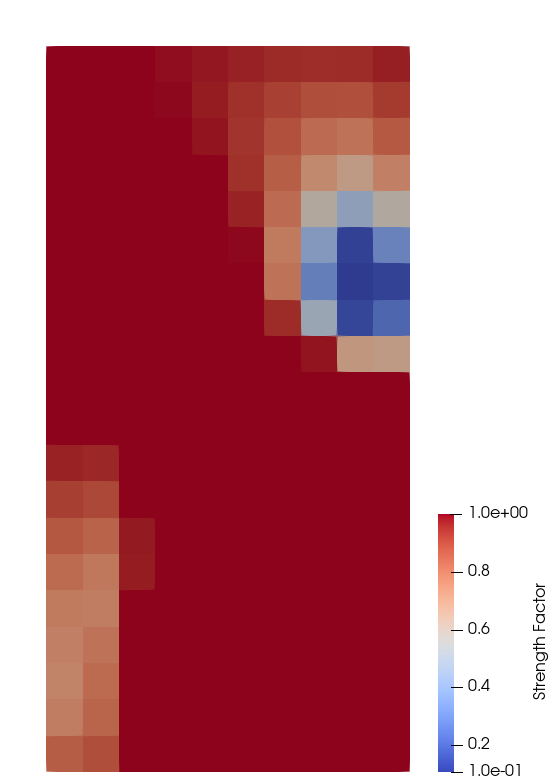}\qquad
	\includegraphics[width=0.25\textwidth]{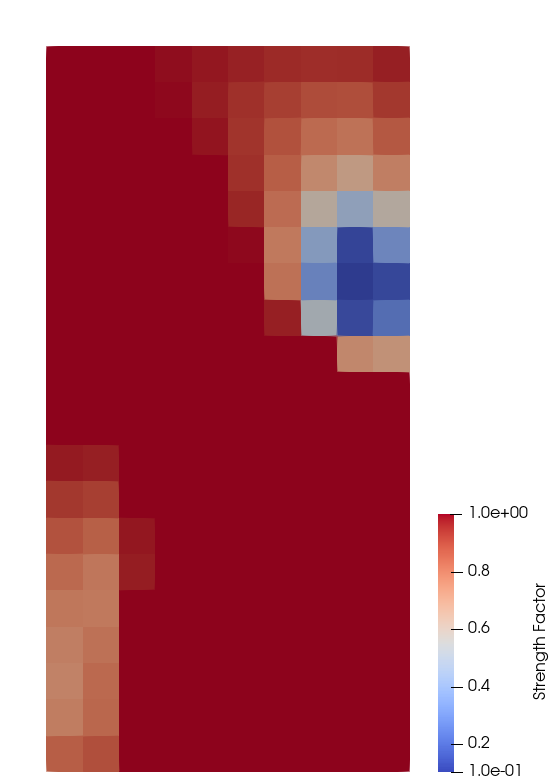}\qquad
	\includegraphics[width=0.25\textwidth]{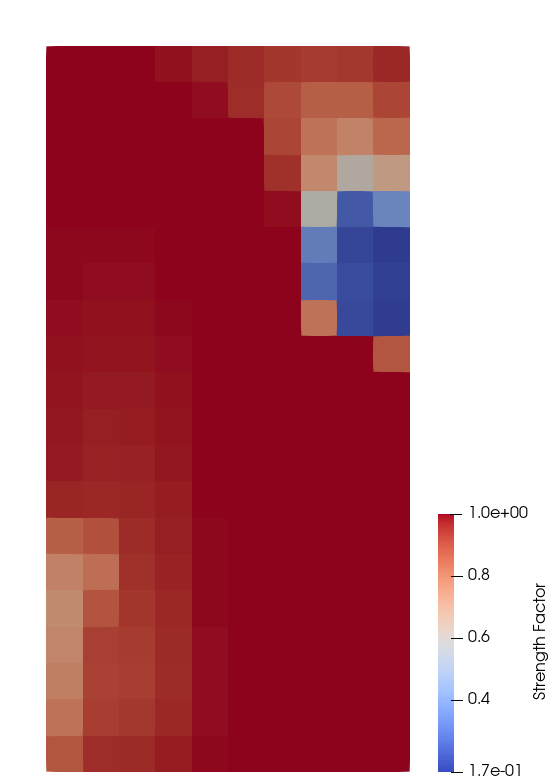}\qquad
	\includegraphics[width=0.25\textwidth]{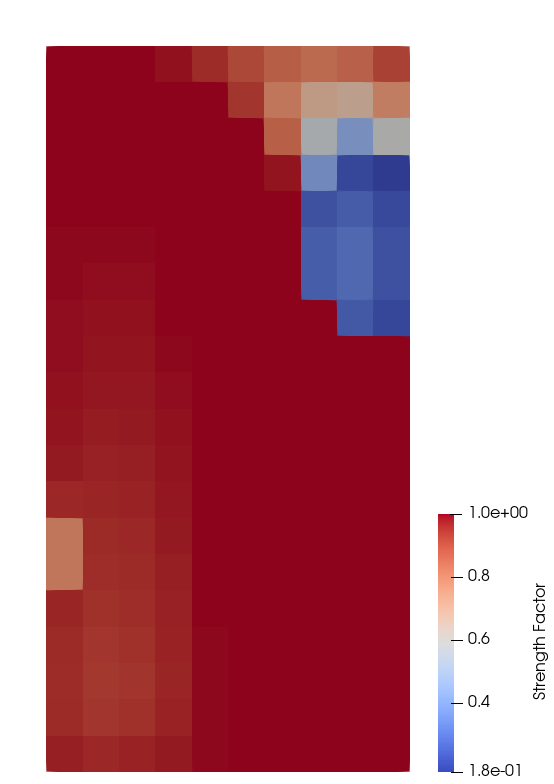}\qquad
    \includegraphics[width=0.25\textwidth]{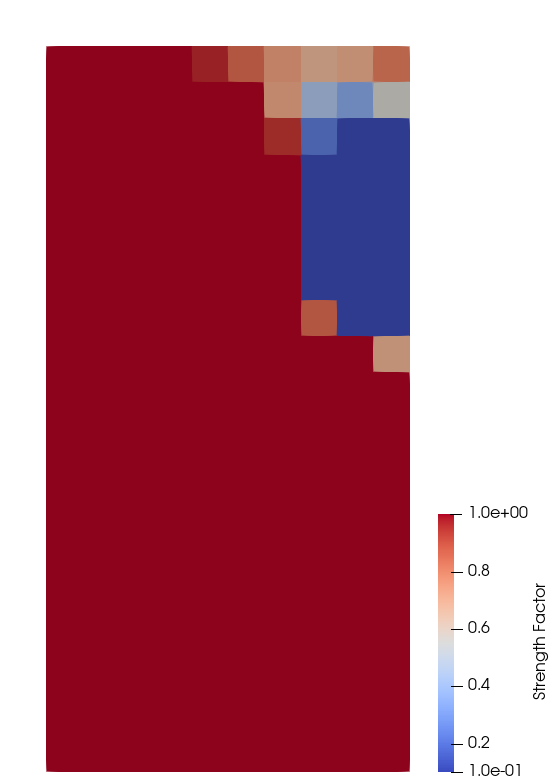}
	\caption{
    Top row: first panel (standard expectation); middle and right panels (CVaR$_\beta$ with
    $\beta = 0.3$ and $0.5$). Bottom row (from left to right): CVaR$_\beta$ with 
    $\beta = 0.6$, $0.7$ and $0.9$. 
    CVaR$_\beta$ performs better than the standard expectation.}
    \label{fig:plate_results2}
\end{figure}

We optimize the same parameters as the previous case, with the difference that our integration limit were now constrained to $\Xi = (0.9,1.1)$. The results are shown in Figure~\ref{fig:plate_results2},  which displays that the risk neutral approach is comparable to  CVaR$_\beta$ for low $\beta$. However, the accuracy improves as $\beta$ gets larger.

\subsection{Example: Large Truss Structure Under Load and Sensor Uncertainty}
\label{s:crane}
\begin{figure}[!hbt]
        \centering
        \includegraphics[width=0.95\textwidth]{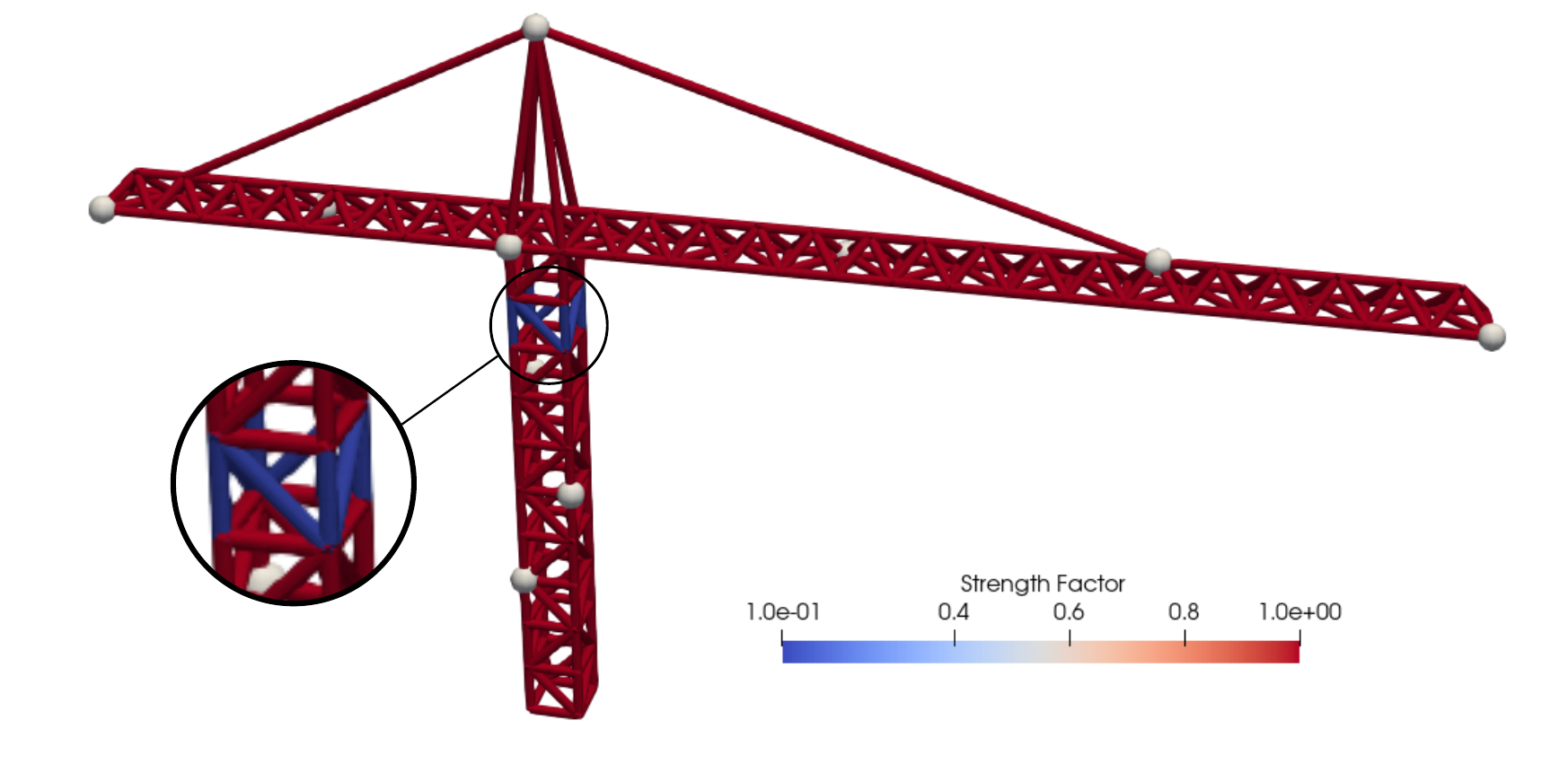}\qquad
        \caption{The goal is to achieve the strength factor shown in this large truss
        configuration. Notice that several beams on the structure have been weakened.}
        \label{f:crane_orig}
\end{figure}

\begin{figure}[!hbt]
  \centering

  \subcaptionbox{Case 1: Downwards load. \\ $\mathbf{-2 \cdot 10^9}$ N on left arm, $\mathbf{-1 \cdot 10^9}$ N on right arm.}%
  {\includegraphics[width=0.49\linewidth]{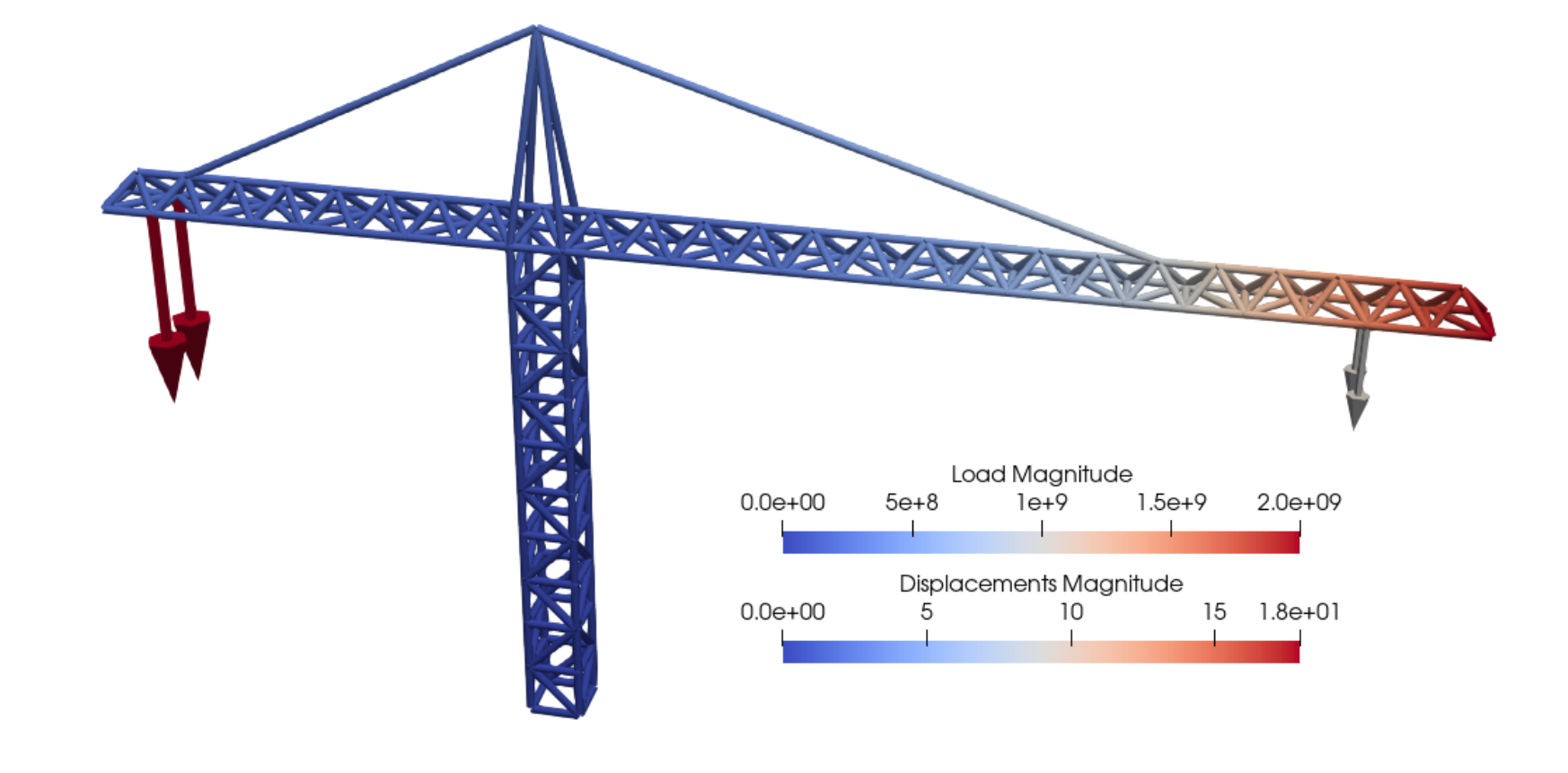}}
  \hspace{\fill}
  \subcaptionbox{Case 2: Side Load. \\ $\mathbf{-2.5 \cdot 10^8}$ N on left arm, $\mathbf{5 \cdot 10^8}$ N on right arm.}%
  {\includegraphics[width=0.49\linewidth]{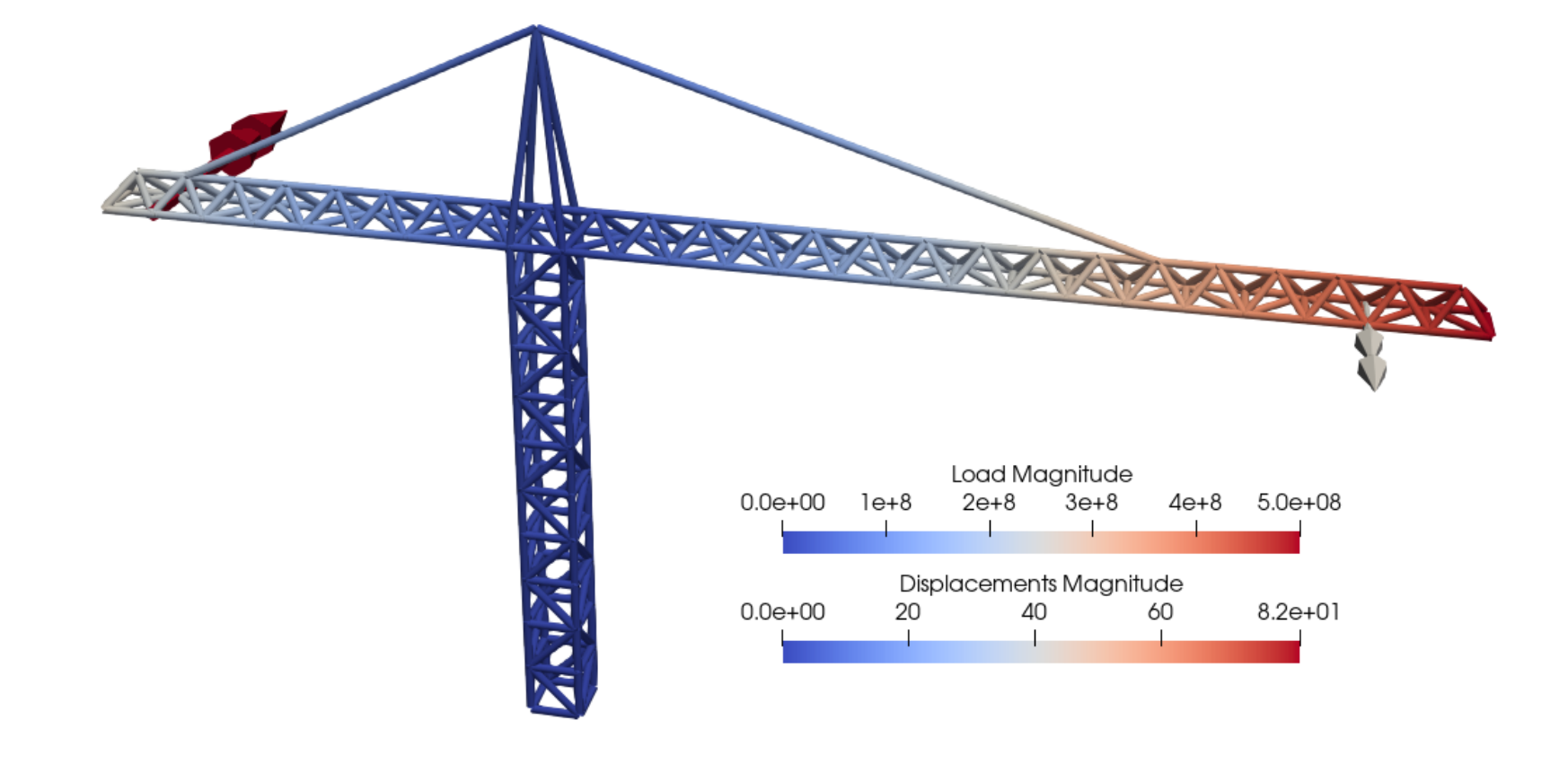}}
  \hspace{\fill}

  \subcaptionbox{Case 3: Downwards load. \\ $\mathbf{-1 \cdot 10^9}$ N on left arm, $\mathbf{-1 \cdot 10^9}$ N on right arm.}%
  {\includegraphics[width=0.49\linewidth]{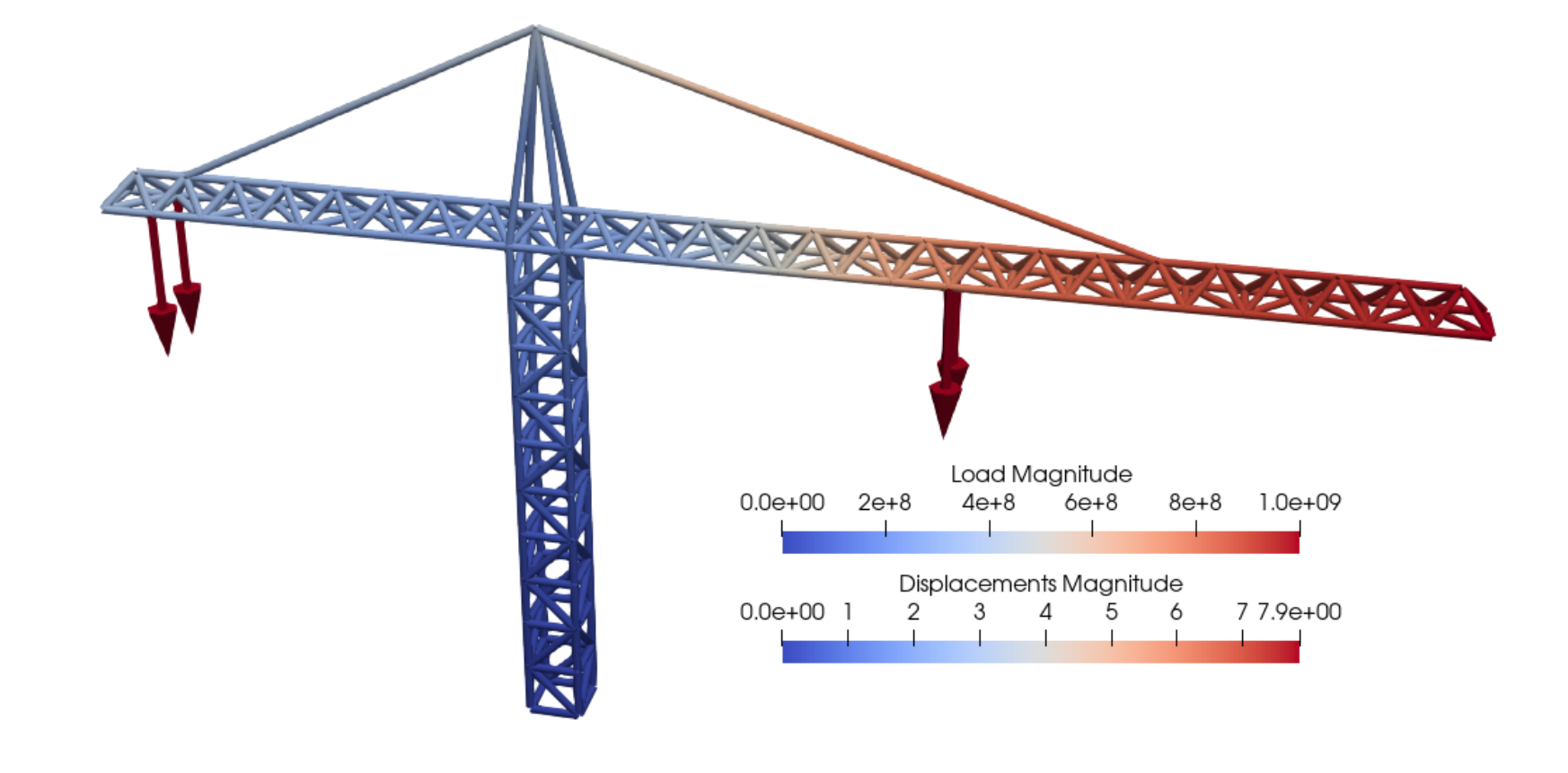}}
  \hspace{\fill}
  \subcaptionbox{Case 4: Side load. \\ $\mathbf{2.5 \cdot 10^8}$ N on left arm, $\mathbf{-5 \cdot 10^8}$ N on right arm.}%
  {\includegraphics[width=0.49\linewidth]{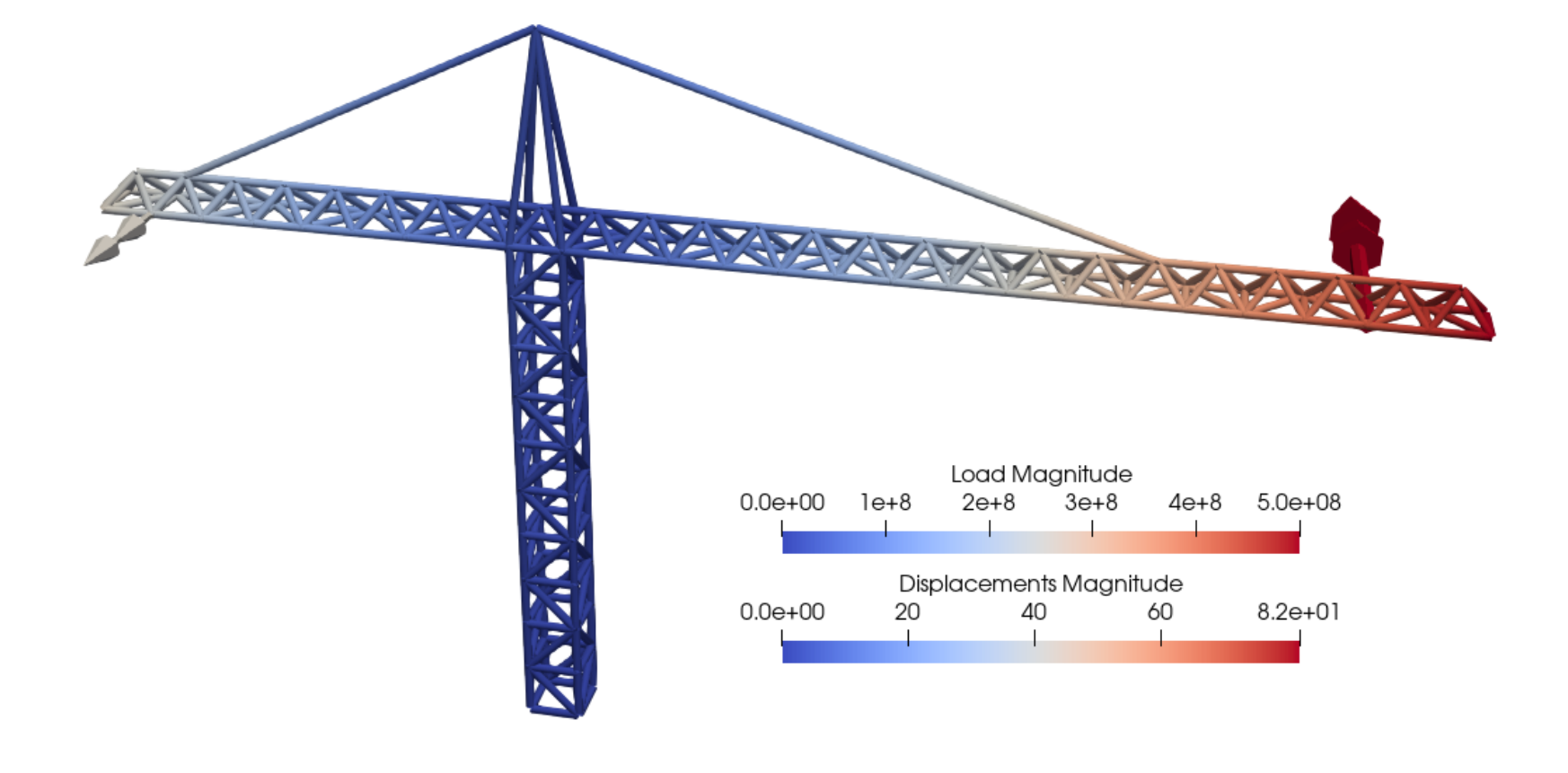}}
  \hspace{\fill}

  \subcaptionbox{Case 5: Torsional load. \\ $\mathbf{\pm 1 \cdot 10^9}$ N on the four points shown.}%
  {\includegraphics[width=0.49\linewidth]{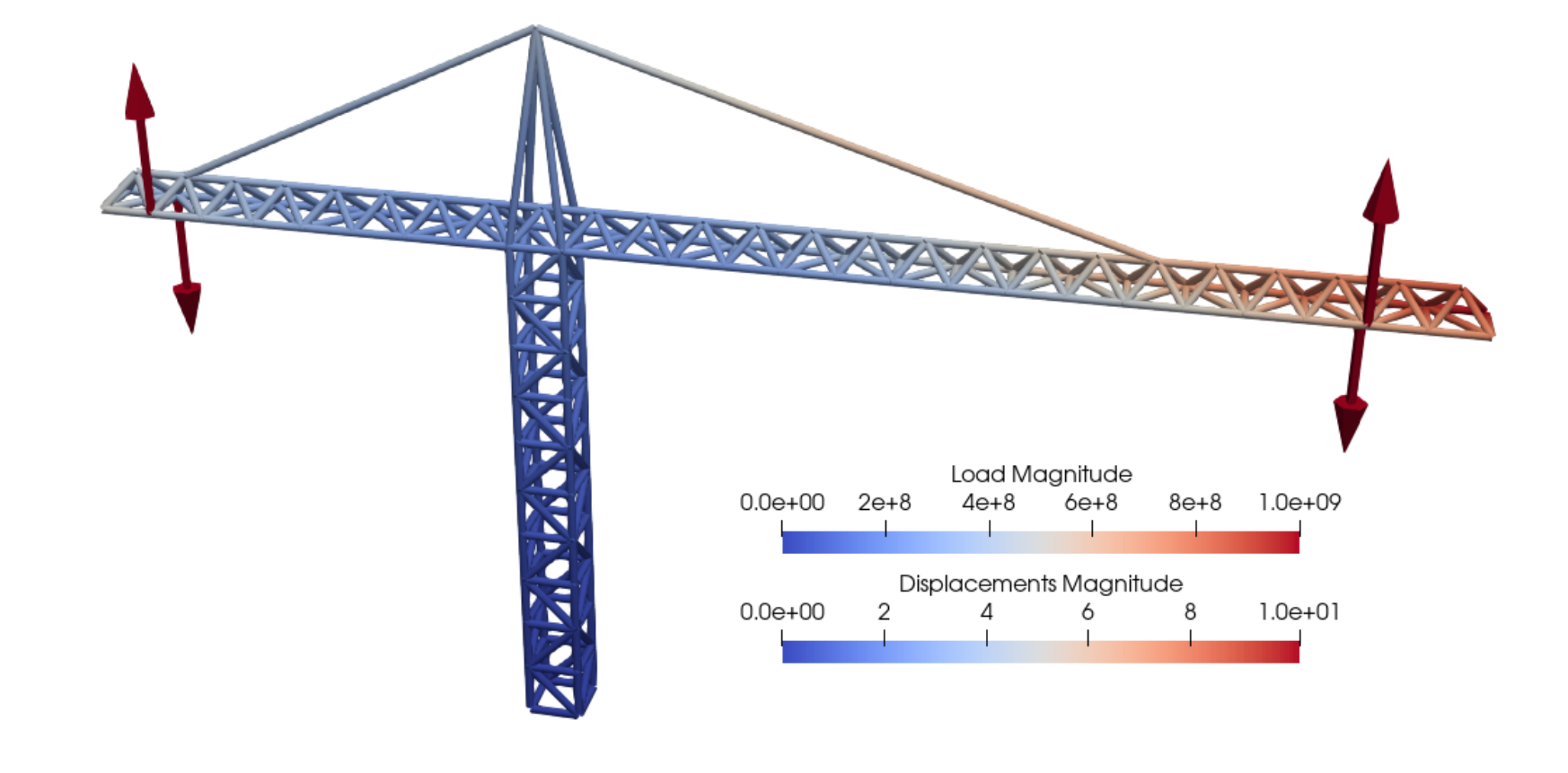}}
  \hspace{\fill}
  \subcaptionbox{Case 6: Torsional load. \\ $\mathbf{\mp 1 \cdot 10^9}$ N on the four points shown.}%
  {\includegraphics[width=0.49\linewidth]{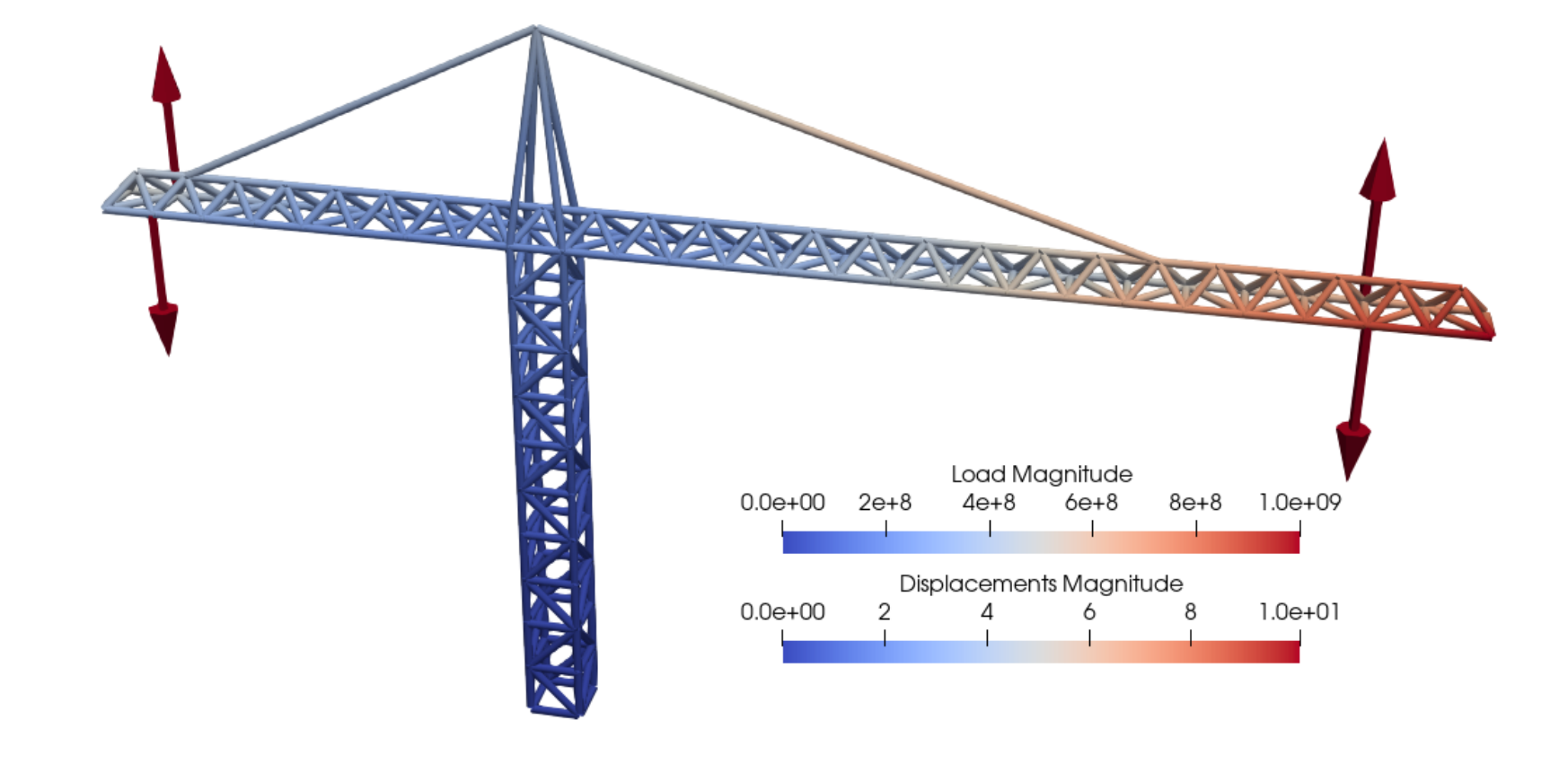}}
  \hspace{\fill}

  \caption{Load cases run simultaneously for the crane problem.}
  \label{fig:crane_cases}
\end{figure}

%
Consider the large truss structure shown in Figure~\ref{f:crane_orig}. Structures such as this one are found
in the solar panels of the international space station,
or in cranes. The structure is composed of 350 beam elements with transversal area $A = 5 \text{cm}^2$. Density, Young’s modulus and Poisson rate were set to $\rho = 7.8$, $E = 2 \cdot 10^{9} $, $\nu = 0.3$, respectively. Moreover, 6~load cases were set up in order to cover a wide set of scenarios. These are presented in Figure \ref{fig:crane_cases}. The loads on each case were allowed to be scaled by a random variable $\xi$, where $\xi$ lies in $\Xi = (0.8,1.2)$.

As in the previous case, the strength factor is discretized using piecewise constant 
elements. Gauss quadrature with 4 terms is used to approximate the integrals 
over $\Xi$. Steepest descent with backtracking line search is used as the 
optimization algorithm. Four smoothing steps are applied to the 
gradient. 

\begin{figure}[!hbt]
        \centering   
	\includegraphics[width=0.49\textwidth]{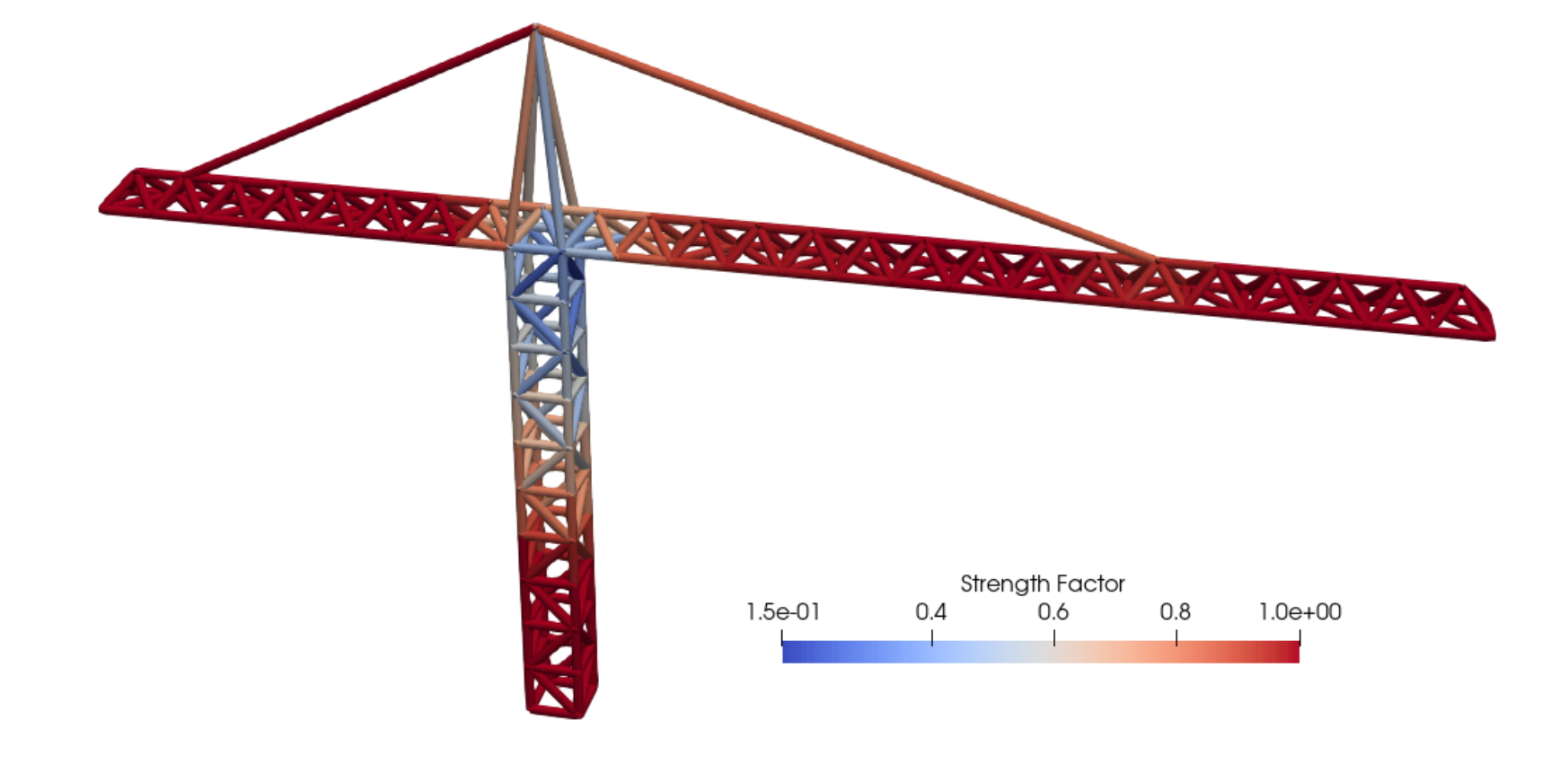} 
	\includegraphics[width=0.49\textwidth]{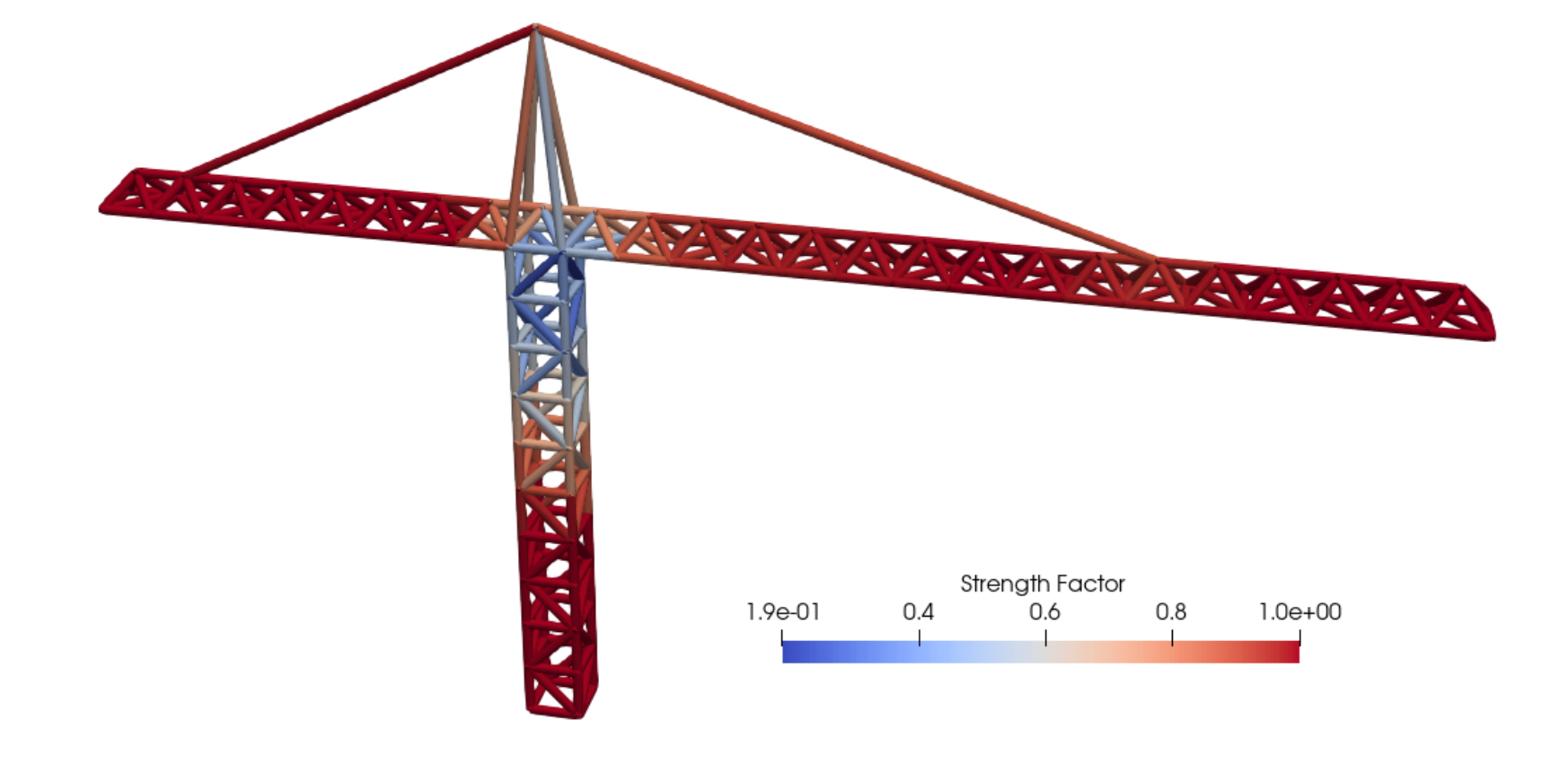}
	\includegraphics[width=0.49\textwidth]{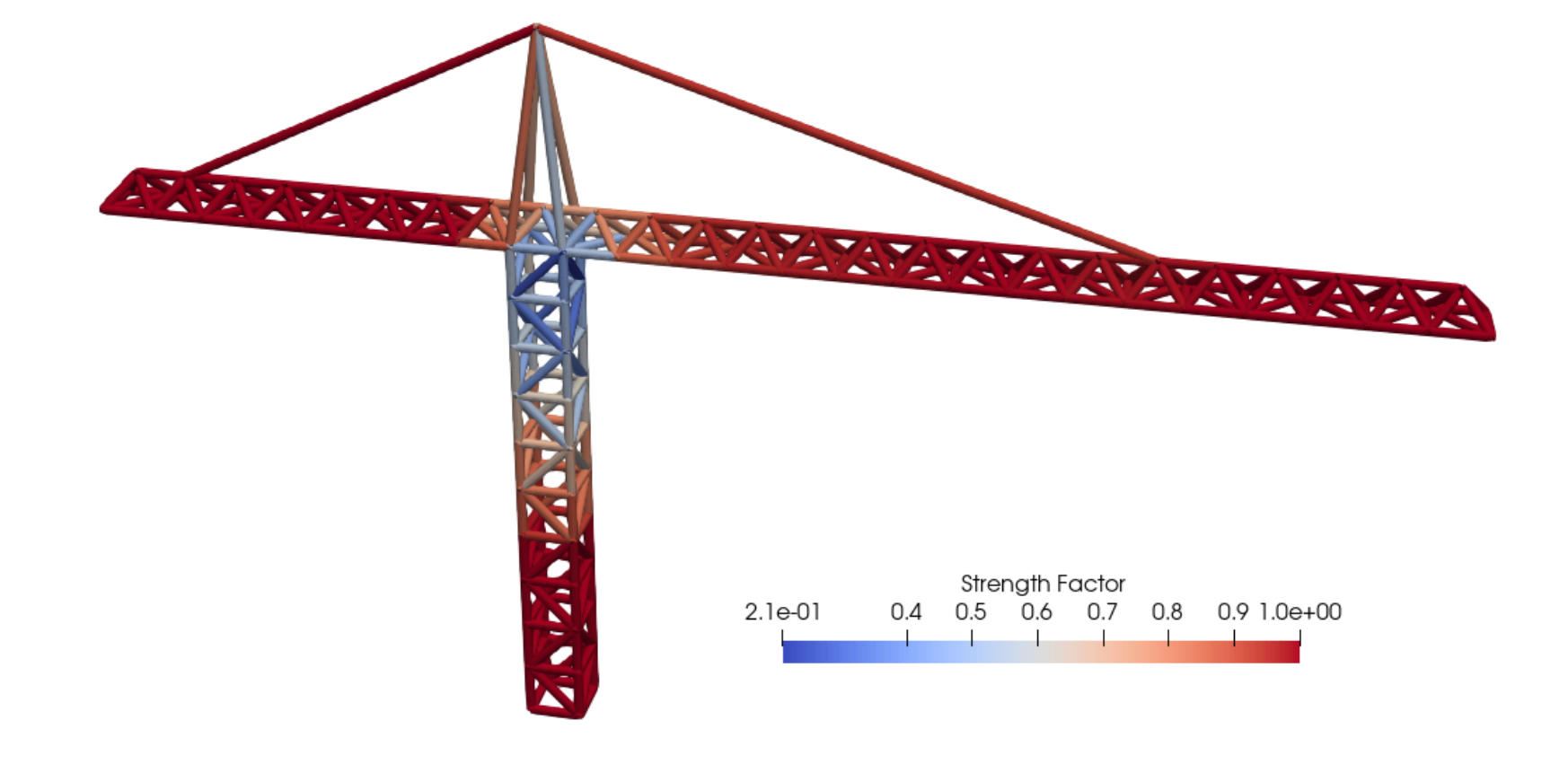}
	\includegraphics[width=0.49\textwidth]{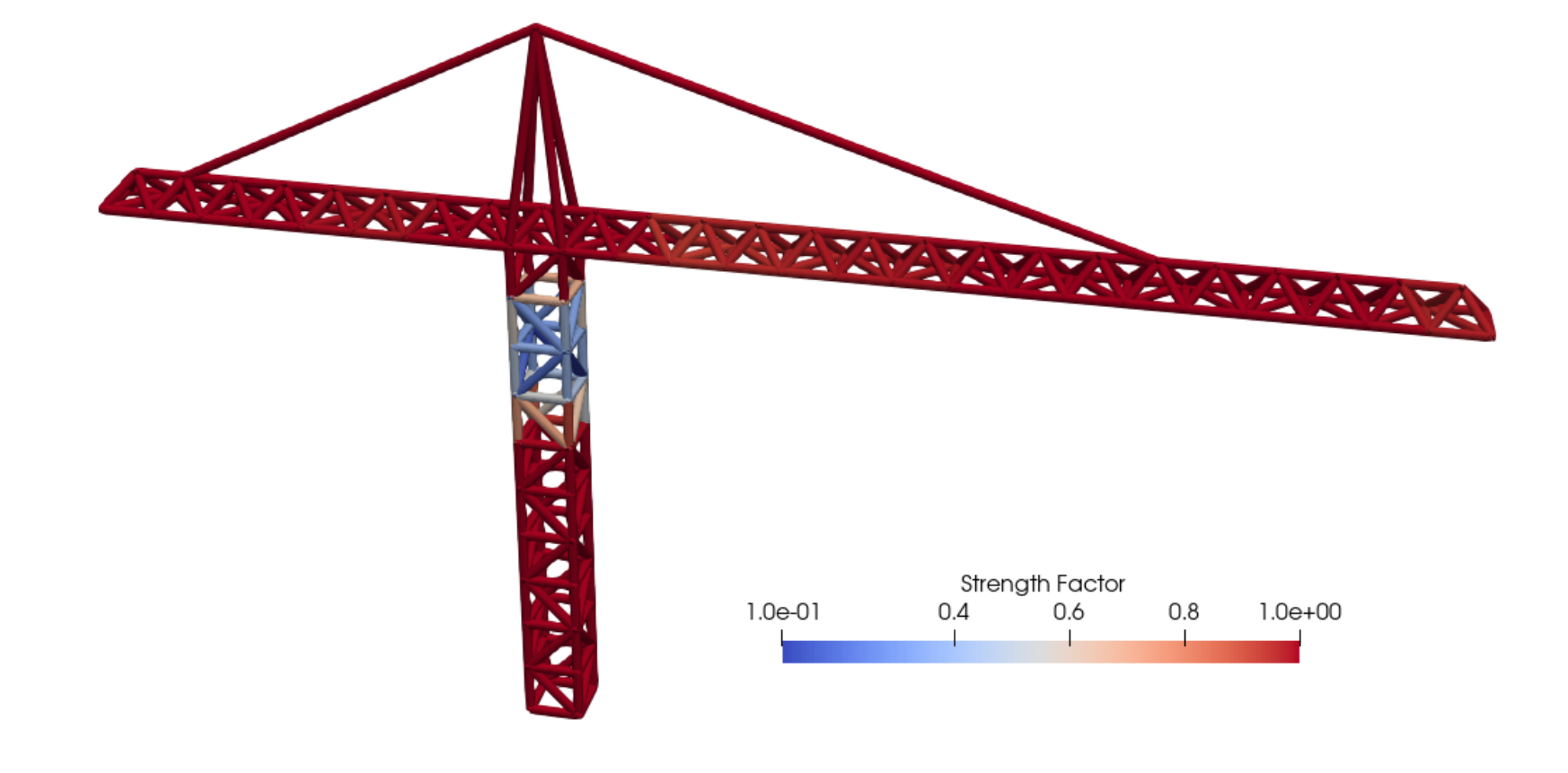}
	\includegraphics[width=0.49\textwidth]{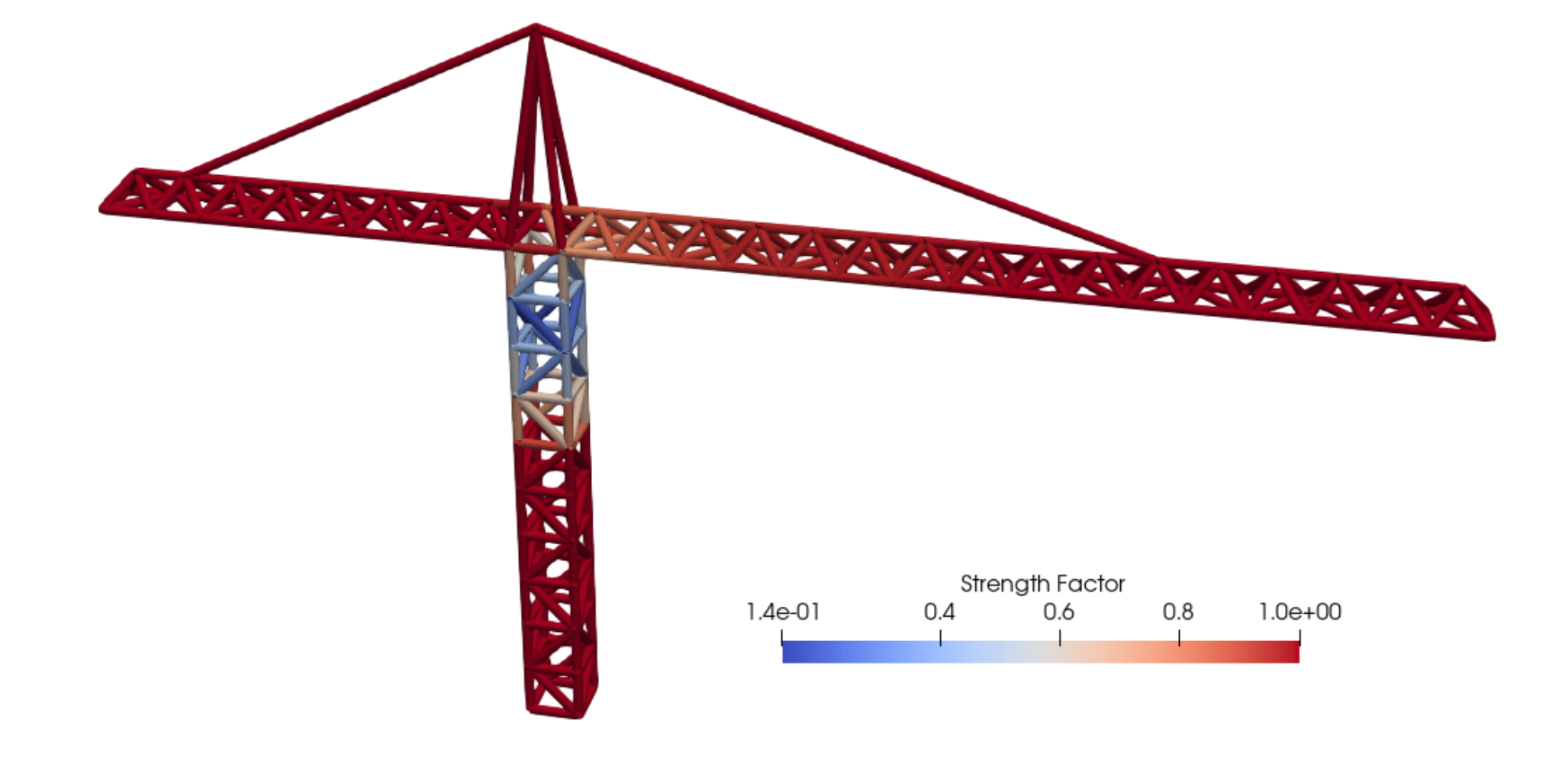}
	\includegraphics[width=0.49\textwidth]{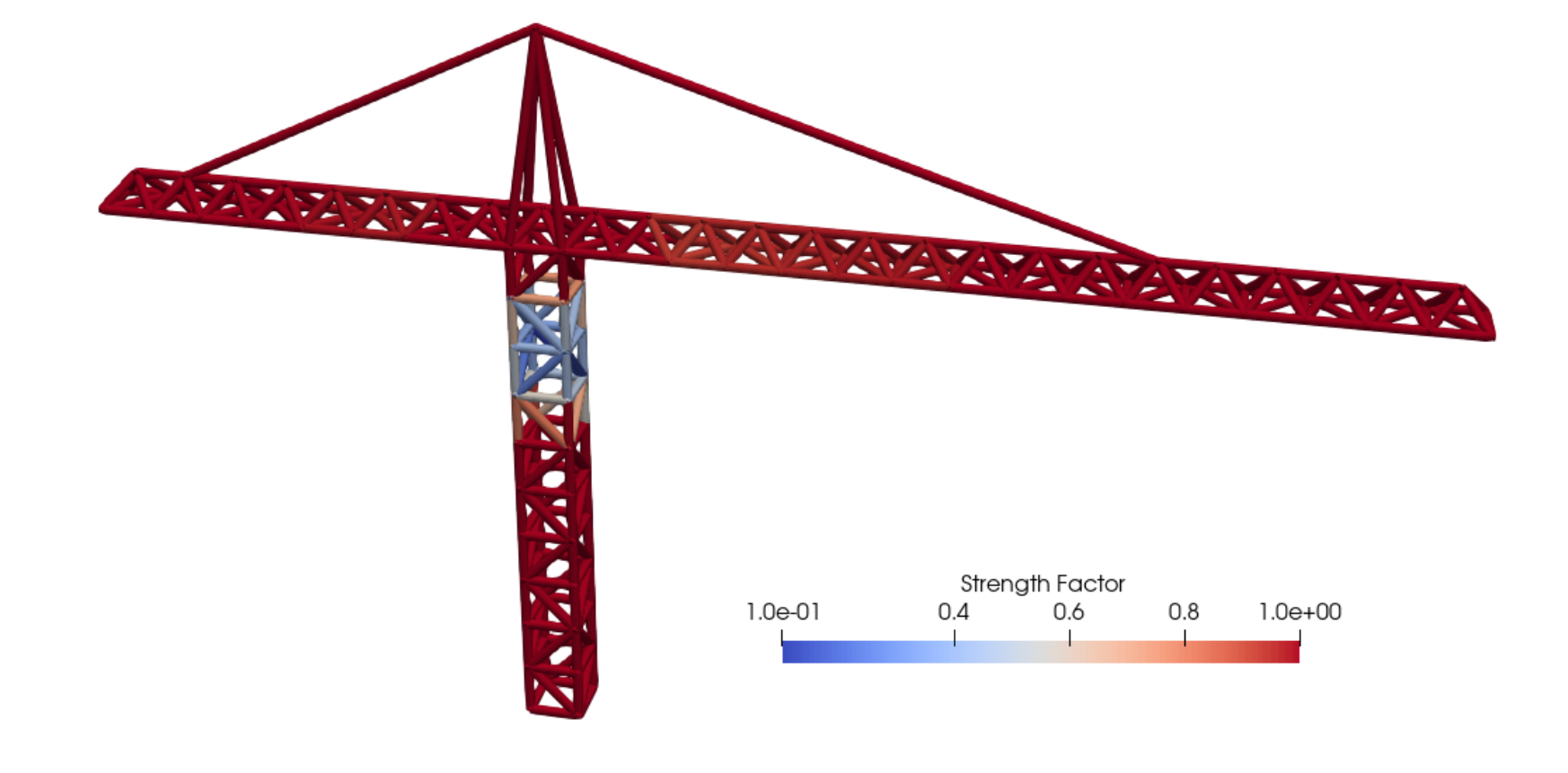}
	\caption{
    Top row: left panel (standard expectation); right panel CVaR$_\beta$ with $\beta = 0.1$. 
    Middle row: CVaR$_\beta$ with $\beta = 0.3$, $0.5$.
    Bottom row: CVaR$_\beta$ with $\beta = 0.7$, $0.8$.
    Clearly, CVaR$_\beta$ performs much better than the standard expectation.}
    \label{f:craneriskaverse}
\end{figure}

For each load case, and each measurement point, a new random variable $\xi$ was generated as in the first plate example. Ten sensors were used in this example (cf. Figure \ref{f:crane_orig}). The weakened configuration was run \texttt{nsensor}*\texttt{ncase} times in order to compute all the target displacements. Figure \ref{f:craneriskaverse} shows a comparison between the standard expectation and CVaR$_\beta$ with different levels of $\beta$. As in the case of the plate, we again observe significant benefits when using CVaR$_\beta$. 

\subsection{Example: Footbridge Under Thermal Loading}
\label{s:footbridge}

\begin{figure}[!hbt]
    \centering
    \includegraphics[width=\textwidth]{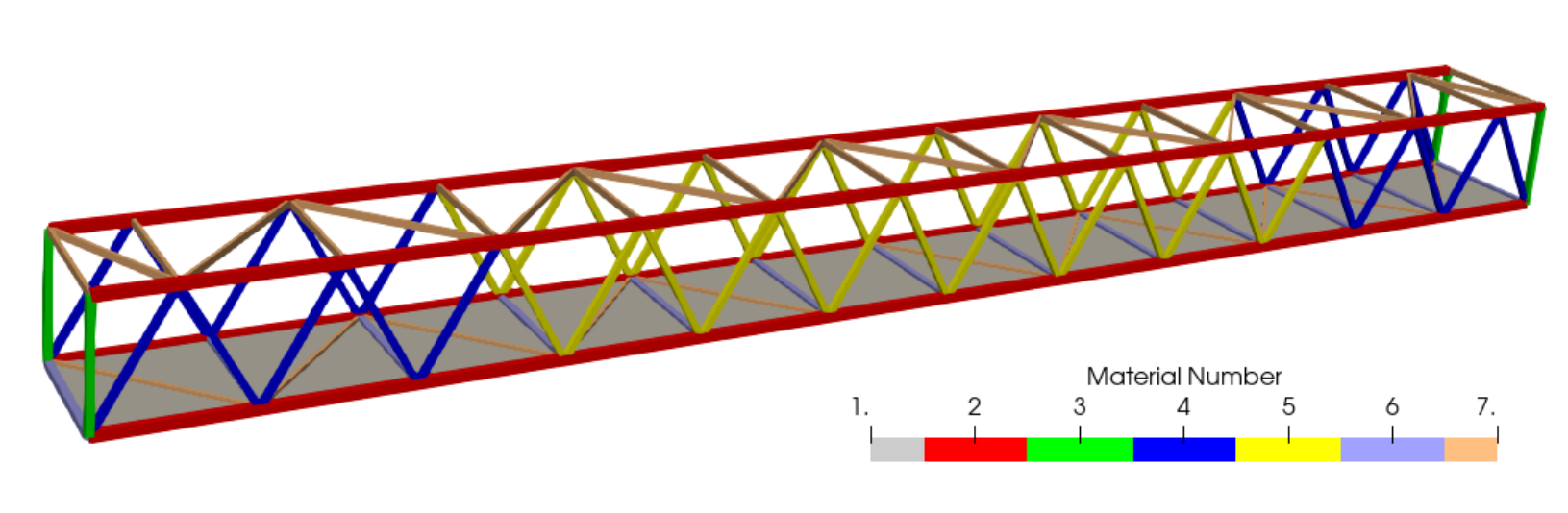}
    \caption{Footbridge: Components.}
    \label{fig:footbridge_components}
\end{figure}

This case considers a typical footbridge and was taken from 
\cite{kilikevivcius2020influence}, see also \cite{FAiraudo_RLoehner_HAntil_2023a} for results in the deterministic setting.
The different types of trusses and plates whose dimensions have been compiled in Table \ref{table:footbridge}, can be 
discerned from Figure \ref{fig:footbridge_components}.
Density, Young's modulus and Poisson ratio were set to
$\rho=7,800~kg/m^3, E=2 \cdot 10^{11}~kg/sec^2/m, \nu=0.3$ respectively.
The structure was modeled using 136 shell and 329 beam elements. 
The bridge is under a distributed load of $1$~MPa in the downwards 
direction, applied to every plate, as well as gravity. 

\begin{table}[!hbt]
\centering
\begin{tabular}{@{}ll@{}}
\toprule
\textbf{Component \#} & \textbf{Shape. Dimensions in mm}                  \\ \midrule
1                    & Steel plate. $t = 10$                                 \\
2                    & Steel beam. Hollow section $300 \times 200 \times 12$ \\
3                    & Steel beam. Hollow section $200 \times 200 \times 10$ \\
4                    & Steel beam. Hollow section $180 \times 180 \times 10$ \\
5                    & Steel beam. Hollow section $180 \times 180 \times 5$  \\
6                    & Steel beam. Hollow section $200 \times 200 \times 10$ \\
7                    & Steel beam. Hollow section $200 \times 100 \times 5$  \\ \bottomrule
\end{tabular}
\caption{Footbridge: Components}
\label{table:footbridge}
\end{table}

In large structures it is typical to have a certain degree of deformation as a result of variations in ambient temperature. This could be a considerable source of uncertainty if this temperature is not correctly accounted for.

In the elasticity equations, a strain component is given by a change in temperature as 
\begin{equation} \label{eq:thermal_strain}
    \epsilon_t = \alpha_{\rm exp} \Delta T,
\end{equation}
where $\alpha_{\rm exp}$ is the coefficient of thermal expansion, which for our case is $\alpha_{\rm exp} = 11 \cdot 10^{-6}$. Naturally, imposing this strain on the constrained system will result in a stress along the structure.

\begin{figure}[!hbt]
    \centering
    \includegraphics[width=\textwidth]{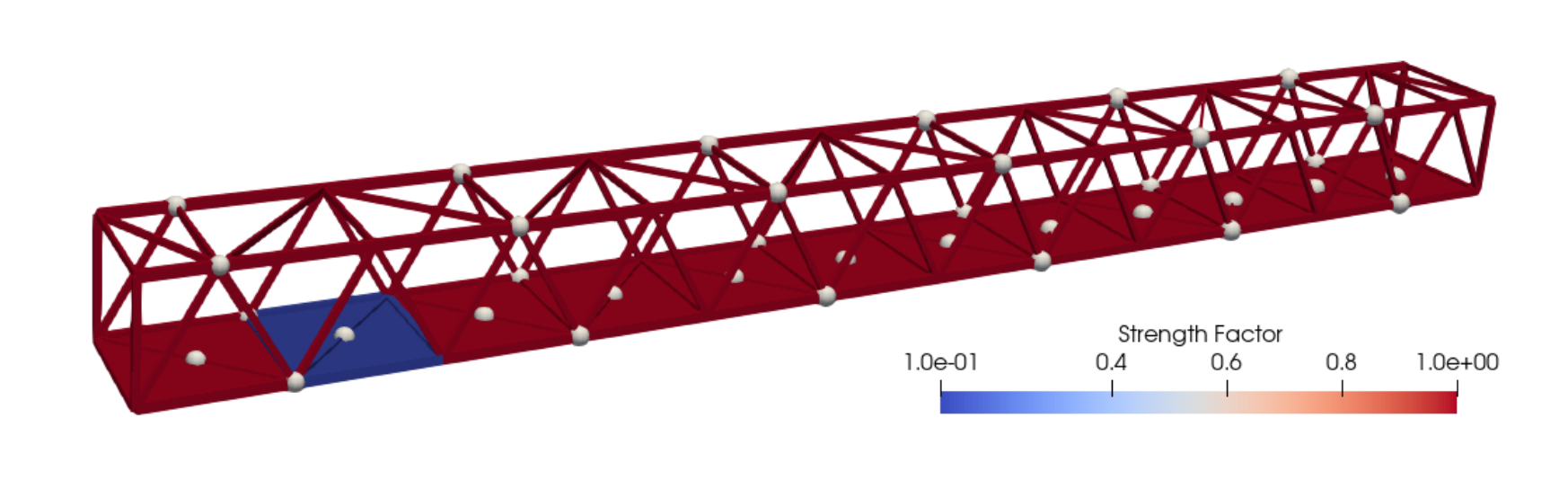}
    \caption{Footbridge: Target configuration with sensors.}
    \label{fig:footbridge_target}
\end{figure}

A weakened configuration such as the one shown in Figure \ref{fig:footbridge_target} is solved for the sensor displacements. In this case however, a $\Delta T = -30K$ is imposed on the whole structure for this target configuration. Gauss quadrature with 4 terms is used to approximate the integrals over $\Xi = (0.9,1.1)$. Steepest descent with backtracking line
search is used as the optimization algorithm. Finally, five smoothing steps are applied to the gradient.

\begin{figure}[!hbt]
    \centering
    \includegraphics[width=\textwidth]{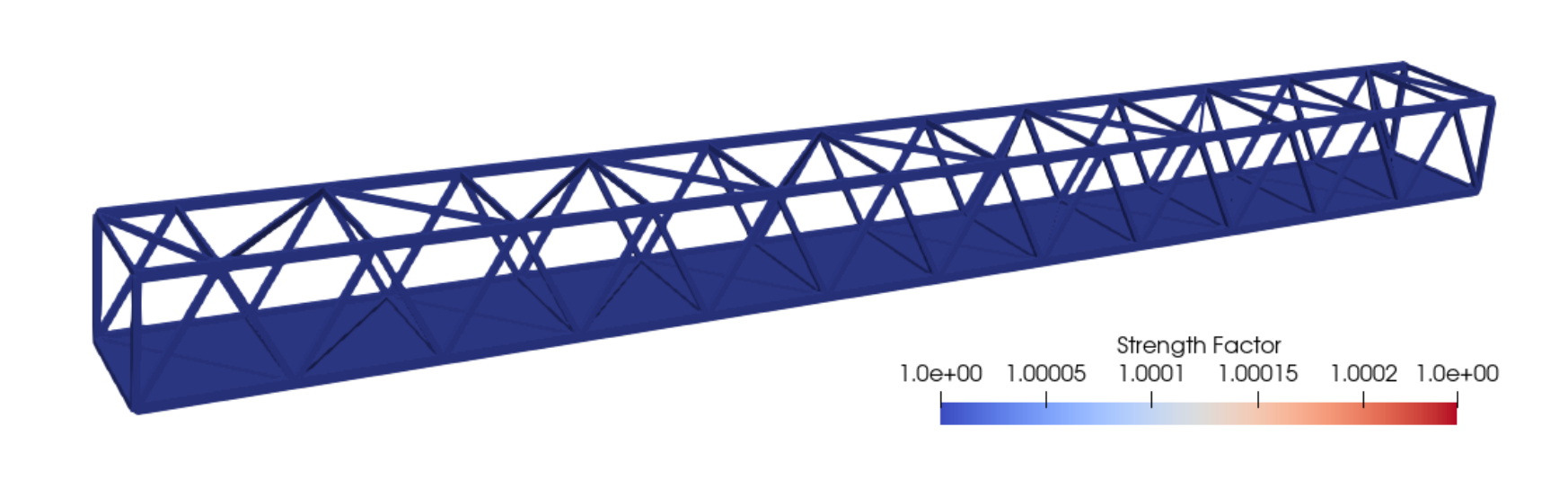}
    \caption{Footbridge: Initial configuration.}
    \label{fig:footbridge_initial}
\end{figure}

We optimize starting with the configuration from Figure \ref{fig:footbridge_initial}, in which $\Delta T = \pm 29.8K$ is assumed to be present as our random variable. This value is $\pm 10\%$ of the ambient temperature $T = 298K$.

\begin{figure}[!hbt]
    \centering
    \includegraphics[width=\textwidth]{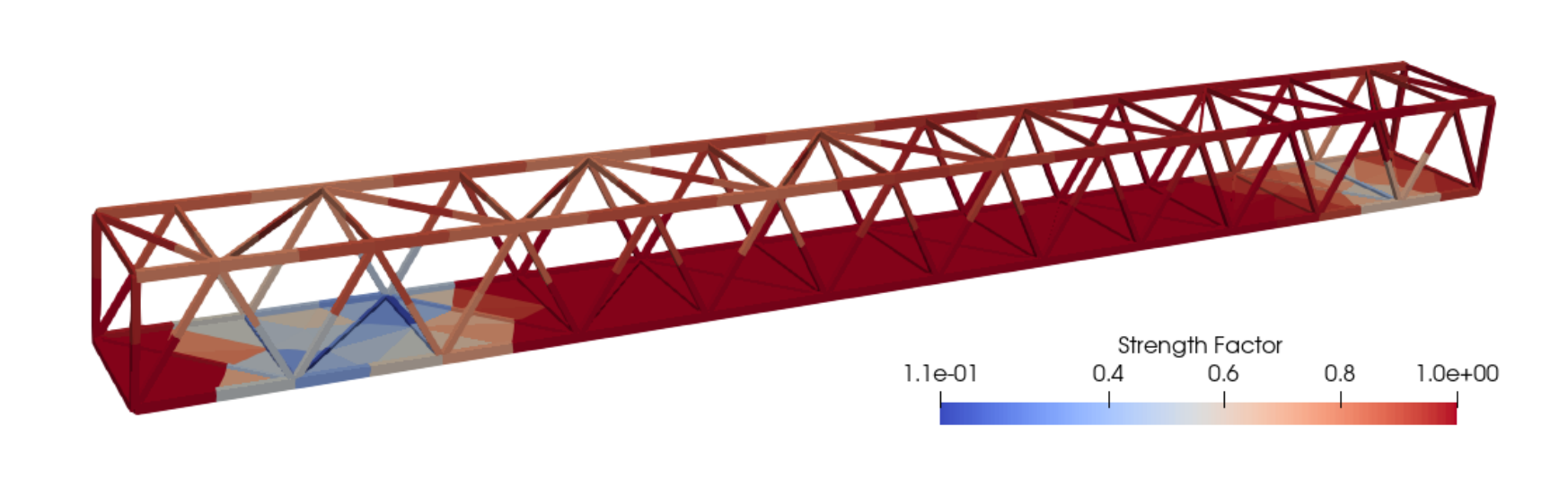}
    \caption{Footbridge: Optimized configuration using risk neutral objective function}
    \label{fig:footbridge_neutral}
\end{figure}

\begin{figure}[!hbt]
    \centering
    \includegraphics[width=\textwidth]{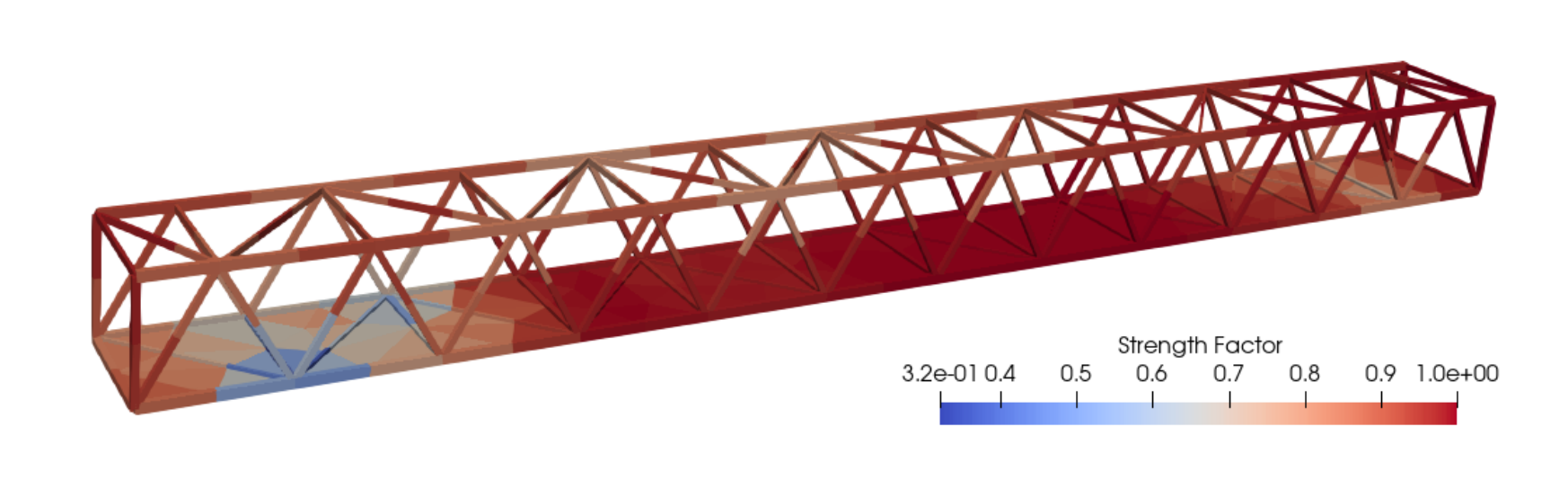}
    \caption{Footbridge: Optimized configuration using CVaR$_\beta$ objective function with $\beta = 0.3$.}
    \label{fig:footbridge_cvar}
\end{figure}

Results are shown for the risk neutral approach in Figure \ref{fig:footbridge_neutral} and for a CVaR$_\beta$ approach with $\beta=0.3$ in Figure \ref{fig:footbridge_cvar}. Both approaches took around 50 gradient descent iterations for a decrease of 3 order of magnitude in the objective.

It is clear that both methods perform well. The uncertain thermal loading results in a \textit{false positive} of a weak spot on the right side of the structure. While the CVaR$_\beta$ method is better at suppressing this spot, it also results in a slightly 
more underdeveloped solution for this case, as can be seen in the color maps.

\section{Conclusion}
\label{s:conclusion}
This work studies the role of risk-measures in digital-twins to tackle
uncertainty in loads and measurements. In particular, the focus is on
identifying weaknesses in structures under uncertainty, leading to an
optimization problem with PDE constraints. This is also termed as `digital twin' because it takes data from the physical system (using sensors) and supplements appropriate predictions. Risk measures such as CVaR$_\beta$ 
are known to be more conservative than risk neutral measures such
as expectation and can lead to robust designs which are resilient 
to rare events.

While the results are very promising, some questions remain unanswered with the use of risk-averse optimization for system identification. As we assume a certain degree of uncertainty on the supposed known loads, and we add more sources of load uncertainty in the shape of thermal or wind loading, the optimization tasks become more complex and costly. The use of a discrete division of the load sources with group loading allows the problems to be tractable, but one could further explore other strategies such as tensor train decomposition \cite{HAntil_SDolgov_AOnwunta_2022b,HAntil_SDolgov_AOnwunta_2023a}. 

As the numerical results show, the conservativeness of the CVaR$_\beta$ approach sometimes leads to underdeveloped solutions, which tends not to be a problem for system identification, as the proper location of the weak spots are found. In addition, it was
found that the best results were obtained for cases in which a large uncertainty in sensor measurements was present. 

\section*{Acknowledgments}

This work is partially supported by NSF grants DMS-2110263, the Air Force Office of Scientific Research (AFOSR) under Award NO: FA9550-22-1-0248, and Office of Naval Research (ONR).

\bibliographystyle{plain}
\bibliography{references.bib}  

\def\ocirc#1{\ifmmode\setbox0=\hbox{$#1$}\dimen0=\ht0 \advance\dimen0
  by1pt\rlap{\hbox to\wd0{\hss\raise\dimen0
  \hbox{\hskip.2em$\scriptscriptstyle\circ$}\hss}}#1\else {\accent"17 #1}\fi}
  \def\cprime{$'$} \def\cprime{$'$}
\begin{thebibliography}{10}

\bibitem{FAiraudo_RLoehner_HAntil_2023a}
Facundo~N. Airaudo, Rainald L{\"o}hner, Roland W{\"u}chner, and Harbir Antil.
\newblock Adjoint-based determination of weaknesses in structures.
\newblock {\em Computer Methods in Applied Mechanics and Engineering},
  417:116471, 2023.

\bibitem{HAntil_LBetz_DWachsmuth_2023a}
Harbir Antil, Livia Betz, and Daniel Wachsmuth.
\newblock Strong stationarity for optimal control problems with non-smooth
  integral equation constraints: Application to a continuous dnn.
\newblock {\em Applied Mathematics \& Optimization}, 88(3):84, 2023.

\bibitem{HAntil_SDolgov_AOnwunta_2022b}
Harbir Antil, Sergey Dolgov, and Akwum Onwunta.
\newblock Ttrisk: Tensor train decomposition algorithm for risk averse
  optimization.
\newblock {\em Numerical Linear Algebra with Applications}, n/a(n/a):e2481.

\bibitem{HAntil_SDolgov_AOnwunta_2023a}
Harbir Antil, Sergey Dolgov, and Akwum Onwunta.
\newblock State-constrained optimization problems under uncertainty: A tensor
  train approach.
\newblock {\em arXiv preprint arXiv:2301.08684}, 2023.

\bibitem{HAntil_DPKouri_MDLacasse_DRidzal_2018a}
Harbir Antil, Drew~P. Kouri, Martin-D. Lacasse, and Denis Ridzal, editors.
\newblock {\em Frontiers in {PDE}-constrained optimization}, volume 163 of {\em
  The IMA Volumes in Mathematics and its Applications}.
\newblock Springer, New York, 2018.
\newblock Papers based on the workshop held at the Institute for Mathematics
  and its Applications, Minneapolis, MN, June 6--10, 2016.

\bibitem{chinesta2020virtual}
Francisco Chinesta, Elias Cueto, Emmanuelle Abisset-Chavanne, Jean~Louis Duval,
  and Fouad~El Khaldi.
\newblock Virtual, digital and hybrid twins: a new paradigm in data-based
  engineering and engineered data.
\newblock {\em Archives of computational methods in engineering}, 27:105--134,
  2020.

\bibitem{dhondt2022calculix}
Guido Dhondt.
\newblock Calculix user’s manual version 2.20.
\newblock {\em Munich, Germany}, 2022.

\bibitem{MHinze_RPinnau_MUlbrich_SUlbrich_2009a}
Michael Hinze, Ren{\'e} Pinnau, Michael Ulbrich, and Stefan Ulbrich.
\newblock {\em Optimization with {PDE} constraints}, volume~23 of {\em
  Mathematical Modelling: Theory and Applications}.
\newblock Springer, New York, 2009.

\bibitem{kilikevivcius2020influence}
Art{\=u}ras Kilikevi{\v{c}}ius, Darius Ba{\v{c}}inskas, Jaroslaw Selech, Jonas
  Matijo{\v{s}}ius, Kristina Kilikevi{\v{c}}ien{\.e}, Darius Vainorius, Dariusz
  Ulbrich, and Dawid Romek.
\newblock The influence of different loads on the footbridge dynamic
  parameters.
\newblock {\em Symmetry}, 12(4):657, 2020.

\bibitem{DPKouri_TMSurowiec_2016a}
D.~P. Kouri and T.~M. Surowiec.
\newblock Risk-averse {PDE}-constrained optimization using the conditional
  value-at-risk.
\newblock {\em SIAM J. Optim.}, 26(1):365--396, 2016.

\bibitem{JLLions_1971a}
Jacques-L. Lions.
\newblock {\em Optimal control of systems governed by partial differential
  equations}.
\newblock Translated from the French by S. K. Mitter. Die Grundlehren der
  mathematischen Wissenschaften, Band 170. Springer-Verlag, New York-Berlin,
  1971.

\bibitem{mainini2015surrogate}
Laura Mainini and Karen Willcox.
\newblock Surrogate modeling approach to support real-time structural
  assessment and decision making.
\newblock {\em AIAA Journal}, 53(6):1612--1626, 2015.

\bibitem{TRRockafellar_SUryasev_2000a}
R~Tyrrell Rockafellar and Stanislav Uryasev.
\newblock Optimization of conditional value-at-risk.
\newblock {\em Journal of risk}, 2:21--42, 2000.

\bibitem{TRRockafellar_SUryasev_2002a}
R~Tyrrell Rockafellar and Stanislav Uryasev.
\newblock Conditional value-at-risk for general loss distributions.
\newblock {\em Journal of banking \& finance}, 26(7):1443--1471, 2002.

\bibitem{AShapiro_DDentcheva_ARuszczynski_2014a}
Alexander Shapiro, Darinka Dentcheva, and Andrzej Ruszczy\'{n}ski.
\newblock {\em Lectures on stochastic programming}, volume~9 of {\em MOS-SIAM
  Series on Optimization}.
\newblock SIAM, Philadelphia, PA; MOS, Philadelphia, PA, second edition, 2014.
\newblock Modeling and theory.

\bibitem{simo2006computational}
Juan~C Simo and Thomas~JR Hughes.
\newblock {\em Computational inelasticity}, volume~7.
\newblock Springer Science \& Business Media, 2006.

\bibitem{FTroeltzsch_2010a}
Fredi Tr{\"o}ltzsch.
\newblock {\em Optimal control of partial differential equations}, volume 112
  of {\em Graduate Studies in Mathematics}.
\newblock American Mathematical Society, Providence, RI, 2010.
\newblock Theory, methods and applications, Translated from the 2005 German
  original by J{\"u}rgen Sprekels.

\bibitem{MUlbrich_2011a}
Michael Ulbrich.
\newblock {\em Semismooth {N}ewton methods for variational inequalities and
  constrained optimization problems in function spaces}, volume~11 of {\em
  MOS-SIAM Series on Optimization}.
\newblock Society for Industrial and Applied Mathematics (SIAM), Philadelphia,
  PA; Mathematical Optimization Society, Philadelphia, PA, 2011.

\bibitem{zienkiewicz2005finite}
Olek~C Zienkiewicz, Robert~Leroy Taylor, and Jian~Z Zhu.
\newblock {\em The finite element method: its basis and fundamentals}.
\newblock Elsevier, 2005.

\end{thebibliography}

\end{document}